\documentclass[a4paper, 10pt]{article}

\usepackage{preamble}
\usepackage{mathtools}

\DeclareMathOperator{\Spec}{Spec}

\title{The Fundamental Theorem of Localizing Invariants}
\author{Victor Saunier}

\begin{document}
\maketitle

\begin{abstract}
    We prove a generalization of the fundamental theorem of algebraic K-theory for Verdier-localizing functors by extending the proof for algebraic K-theory of spaces to the realm of stable $\infty$-categories. The formula behaves much better for Karoubi-localizing functors, the Verdier-localizing invariants which are additionally invariant under idempotent completion. \\
    
    This general fundamental theorem specializes to new formulas in the context of non-connective K-theory, topological Hochschild homology and topological cyclic homology as well as connective K-theory of stable $\infty$-categories, and generalizes several known formulas for algebraic K-theory of spaces or connective K-theory of ordinary rings, ring spectra, schemes and $\S$-algebras.
\end{abstract}

\vspace{2em}

\tableofcontents

\newpage
\section{Introduction}
\vspace{0.5em}

\hspace{1.2em} In algebraic K-theory of rings, the celebrated fundamental theorem, also known as the Bass-Heller-Swan formula, states that, for any ring $R$ and any integer $n\in\Z$, there is a natural isomorphism:
$$
    K_n(R[t, t^{-1}]) \cong K_n(R)\oplus K_{n-1}(R)\oplus (N_+K)_n(R)\oplus (N_-K)_n(R)
$$
where $N_\pm K(R)$ are so called \textit{nilterms}, isomorphic functors obtained as the kernel of the split morphism $K_n(R)\to K_n(R[t])$ which vanish when $R$ is a regular ring. The $K_{-1}(R)$ appearing in degree zero is the first group of negative K-theory of rings as defined by Bass, and one can think of this formula as providing a iterative definition for the negative K-groups (this was actually the historical appraoch by Bass). This result was proved by Bass-Heller-Swan for $K_0$, $K_1$ and negative K-groups and later in full generality for non-negative integers by Quillen, whose proof was presented by Grayson in \cite{Grayson}. It was then extended to schemes by Thomason-Trobaugh in \cite{ThomasonTrobaugh}. \\

In \cite{HKVWW}, Hüttemann, Klein, Vogell, Waldhausen and Williams proved a similar equivalence for the finitely-dominated variant of $A(X)$ the algebraic K-theory of a space $X$, originally defined by the fourth author. If $X$ is a space, their version of A-theory, $A^{fd}(X)$, is defined to be the K-theory of the Waldhausen category $\Rcal_{fd}(X)$ of finitely-dominated retracts of $X$. They showed in \textit{loc. cit.} the following splitting:
$$
    A^{fd}(X\times S^1)\simeq A^{fd}(X)\times\Bcal A^{fd}(X)\times N_+A^{fd}(X)\times N_-A^{fd}(X)
$$
where $N_\pm A^{fd}(X)$ are homeomorphic nilterms and $\Bcal A^{fd}(X)$ is a non-connective but canonical delooping of $A^{fd}(X)$. The nilterm are also obtained as some kernel, though there is no space that plays the role of $R[t]$ in this context. This can be thought as an extension of the fundamental theorem to some class of ring spectra\footnote{By which we mean $E_{1}$-ring spectra}: indeed, $A(X)$ is equivalently defined as the K-theory of the ring spectrum $\S[\Omega X]$ and $\S[\Omega (X\times S^1)]\simeq\S[\Omega X][t, t^{-1}]$. The formula above is then simply the fundamental theorem for these "brave new rings". \\

Other extensions of this fundamental theorem have been made in recent years, expanding the known cases: Lück and Steimle have proven a twisted formula for both connective and non-connective K-theory of additive categories in \cite{LuckSteimle} and Fontes and Ogle have shown the theorem holds for connective $\S$-algebras in \cite{FontesOgle}. Moreover, Hüttemann has extended the result for strongly $\Z$-graded rings in \cite{Huttemann}. 

More recently, the result has been shown to hold for a greater class of invariants than just flavors of K-theory. In \cite{Tabuada}, Tabuada proved the formula held with vanishing nilterms for every $\A^1$-homotopy invariant, derived Morita-invariant functor from dg-categories which was furthermore localizing. In the more geometric context of spectral algebraic geometry, Cisinski and Khan showed that for localizing invariants of stable $R$-linear $\infty$-category, where $R$ is a connective ring spectra, the formula held as well (see Theorem 4.3.1 in \cite{CisinskiKhan}). \\

The motivation of this article is to provide a general statement which encompasses many of the above formulas, similarly to the last two cited papers. Lured by recent development, such as \cite{Blumberg} or \cite{Barwick}, which have made clear the practicality and the usefulness of higher category theory for K-theoretic purposes, we will prove a higher categorical version of the fundamental theorem and, as in \cite{CisinskiKhan}, show that it not only holds for algebraic K-theory but also for a whole array of invariants which satisfy a key property of \textit{localization}, the precise flavor of it we will explain in the following. Amongst such localizing invariants are notably topological Hochschild homology, topological cyclic homology and non-connective K-theory. In particular, our main result recovers as a special case the fundamental theorem of algebraic K-theory for schemes as in \cite{ThomasonTrobaugh} and connective $\S$-algebras as in \cite{FontesOgle}, the A-theoretic version of \cite{HKVWW} and provides a different proof of the central case of \cite{CisinskiKhan}. We do not however generalize the twisted versions of \cite{LuckSteimle} and \cite{Huttemann}.

In contrast to \cite{CisinskiKhan} however, our proof avoids references to $\E_\infty$-ring spectra or more technical ideas of spectral algebraic geometry (in particular, it does not rely on any result of \cite{SAG}) and is purely contained in the theory of stable $\infty$-categories. This allows for finer control and thus slightly refined statements, which involve stable $\infty$-categories which may not be idempotent complete, unlike Cisinski and Khan. \\

Finally, let us mention that in the upcoming fourth part \cite{HermKIV} of the series of hermitian K-theory (\cite{HermKI, HermKII, HermKIII}), the authors prove a similar Bass-Heller-Swan formula for hermitian K-theory and more generally for Poincaré-Verdier localizing functors. However, due to technicalities arising from the setting of Poincaré $\infty$-categories, the proof they supply differs from the classical proof, including the one we present here. Notably, it does not feature a projective line. \\

\textbf{Notations and conventions.} As we explained, the higher categorical setting is one of the motivation of this article. Thus we adopt throughout this paper the language of $\infty$-categories as developed by Jacob Lurie in \textit{Higher Topos Theory} \cite{HTT} and \textit{Higher Algebra} \cite{HA}. We recall now the main concepts and notations we will use. \\

A stable $\infty$-category is an $\infty$-category with a zero object — we say such $\infty$-categories are \textit{pointed} — such that every morphism has a fiber and a cofiber and additionally, that fiber and cofiber sequences coincide. Stable $\infty$-categories admit all finite limits and colimits but also, cartesian squares coincide with cocartesian squares — and are called exact squares. Functors preserving finite limits and colimits are called exact functors. \\

We will denote $\cat{Cat}_\infty$ the $\infty$-category of (small) $\infty$-categories, and $\CatEx$ the non-full subcategory of stable $\infty$-categories and exact functors between them. When $\Ccal$ is an $\infty$-category, we denote $\Ind(\Ccal)$ its $\Ind$-construction, $\Ccal^c$ its full subcategory of compact object and $\Idem(\Ccal)$ its idempotent completion. In particular, \cite{HTT} 5.4.2.4 gives an equivalence $\Ind(\Ccal)^c\simeq\Idem(\Ccal)$. When $\Ccal$ is stable, so are all of the above $\infty$-categories. \\

We let $\CatExIdem$ denote the full subcategory of $\CatEx$ of idempotent-complete stable $\infty$-category. Finally, $S$ will denote the $\infty$-category of spaces and $\Sp$ the stable $\infty$-category of spectra, which is the stabilisation of the former. \\

\textbf{Main results.} In the original proof of \cite{Grayson}, the fundamental theorem for algebraic K-theory of rings is deduced from a property of K-theory regarding localisation of rings. We proceed by a categorification of this notion, following the ideas of \cite{Blumberg}, which have realized (non-connective) algebraic K-theory as a functor built from the universal \textit{localizing invariant}.

To properly capture the exact flavour of this localization, we adopt however a different semantic from that of \cite{Blumberg} or \cite{CisinskiKhan}, which is inspired from the series of papers on hermitian K-theory for higher categories \cite{HermKI, HermKII, HermKIII, HermKIV}. First, let us define three classes of cofibers of $\CatEx$:

\begin{defi}
    A \textit{Verdier sequence} is a sequence of $\CatEx$
    $$
        \begin{tikzcd}
            \Ccal\arrow[r, "i"] & \Dcal\arrow[r, "p"] &\Ecal
        \end{tikzcd}
    $$
    which is both a fiber and a cofiber in $\CatEx$. Exact functors $i$ fitting in such sequences are called Verdier inclusions and $p$ Verdier projections. \\
    
    If $p$ has a left adjoint (resp. right), then the Verdier sequence is called left-split (resp. right-split). Sequences that are both left- and right-split are called \textit{split-Verdier sequences}. \\
    
    A \textit{Karoubi sequence} is a sequence of $\CatEx$
    $$
        \begin{tikzcd}
            \Ccal\arrow[r, "i"] & \Dcal\arrow[r, "p"] &\Ecal
        \end{tikzcd}
    $$
    which is sent by the idempotent-completion functor $\Idem$ to both a fiber and a cofiber in $\CatExIdem$, the full subcategory of $\CatEx$ of stable idempotent-complete $\infty$-categories. Exact functors $i$ fitting in such sequences are called Karoubi inclusions and $p$ Karoubi projections. 
\end{defi}

We will characterize Verdier projections and Karoubi projections; in particular, we will see that Verdier sequences are in particular Karoubi sequences. Functors $\CatEx\to\Sp$ sending either of those three classes of cofibers to exact sequences of spectra will be our different flavours of localizing invariants:

\begin{defi}
    A reduced functor $F:\CatEx\to\Sp$ is called:
    \begin{itemize}[itemsep=0em]
        \item \textit{additive} or \textit{split-Verdier localizing} if it sends split-Verdier sequences to exact sequences
        \item \textit{Verdier localizing} if it sends Verdier sequences to exact sequences
        \item \textit{Karoubi localizing} if it sends Karoubi sequences to exact sequences
    \end{itemize}
\end{defi}

By the remark above, Karoubi-localizing functors are in particular Verdier-localizing. In fact, we will see that they are exactly the Verdier-localizing functors that are invariant under idempotent completion in Proposition \ref{Karoubi=Verdier+IdempotentInvariance}.

Note that we do not ask that localizing invariants preserve filtered colimits. Although this is the case of many of the examples we will encounter, such as algebraic K-theory or topological Hochschild homology, it will not be needed to prove our results. In this, we differ from the definitions of \cite{Blumberg} and adopt the conventions of \cite{HermKII} that separate the two properties.

Our definitions being only slightly different from theirs, the full results of \cite{Blumberg} could be tweaked to apply to our setting. However, for simplicity and because they do not play a part in the further parts of this article, we have chosen to omit the discussion of non-commutative motives and associated results. \\

In our categorified setting, the ring of Laurent polynomials $R[t, t^{-1}]$ appearing in the Bass-Heller-Swan formula will be replaced by some stable $\infty$-category that we preemptively denote $S^1\otimes\Ccal$. Beware that what we denote $\otimes$ is not the usual tensor product of symmetrical monoidal structures on $\cat{Cat}_\infty^{Ex}$ or $\cat{Cat}_\infty$ but a hybrid version of the two. \\

First recall the usual tensor product of $\CatEx$. If $\Ccal$, $\Dcal$ are stable $\infty$-categories, then there exists a stable $\infty$-category which is universal for functors from $\Ccal\times\Dcal$ which are exact in both variables. Since we are interested in tensoring a stable $\Ccal$ by any $\infty$-category $K$, which need not be stable, we will drop the exactness requirement on one side. Hence, our $K\otimes\Ccal$ is universal for functors from $K\times\Ccal$ which are exact only in the second variable. 

Notice that when $\Ccal$ and $\Dcal$ are idempotent complete, their tensor product need not be and the internal tensor product of $\CatExIdem$ is obtained by taking the idempotent completion of the former. Similarly, when $\Ccal$ is idempotent complete, the tensoring $K\otimes\Ccal$ need not be but we can consider its idempotent completion $K\hat\otimes\Ccal$. Let $R$ be a ring spectrum and denote $\Perf(R)$ the stable $\infty$-category of compact $R$-modules, then we have an equivalence $S^1\hat\otimes\Perf(R)\simeq\Perf(R[t, t^{-1}])$. Indeed, this follows from the explicit construction of Proposition \ref{TensorConstruction}, which realizes $S^1\hat\otimes\Perf(R)$ as the compact objects of $\Fun((S^1)^{\mathrm{op}}, \Mod_R)$, and the identification of $R[t, t^{-1}]$-modules as $R$-modules with a $\Z$-action. \\

Denote $S^1_\pm$ the category $B\N$, where the $\pm$ comes from either identifying $\N$ as non-negative or non-positive integers in $\Z$, and let $N_\pm F(\Ccal)$ be the cokernel of $F(\Ccal)\to F(S^1_\pm\otimes\Ccal)$ for a stable $\Ccal$. Our first main result relates $F(S^1\otimes\Ccal)$, $N_\pm F(\Ccal)$ and $F(\Ccal)$ for $F$ Verdier-localizing and $\Ccal$ idempotent complete. 

\begin{thm} \label{MainResultIntro}
    Let $\Ccal$ be a stable idempotent-complete $\infty$-category and $F:\CatEx\to\Sp$ a Verdier-localizing invariant, then, we have the following equivalence of spectra:
    $$
        F(S^1\otimes\Ccal)\simeq F(\Ccal)\oplus\Sigma F(\Ccal)\oplus N_+F(\Ccal)\oplus N_-F(\Ccal)
    $$
    where $N_\pm F(\Ccal)$ are the nil-terms previously defined.
\end{thm}

The above formula is not quite the Bass-Heller-Swan formula, even for rings: even though we supposed $\Ccal$ idempotent complete, $S^1\otimes\Ccal$ need not be; such is the case of $\Ccal=\Perf(R)$. To get the correct statement, one needs to replaced $S^1\otimes\Ccal$ by its idempotent-completion $S^1\hat\otimes\Ccal$. This is automatic in the special case of Karoubi-localizing invariants, which is exactly further supposing invariance under idempotent completion.

\begin{thm}[Generalized Bass-Heller-Swan formula] \label{MainResult2Intro}
    Let $\Ccal$ be a stable $\infty$-category and $F:\CatEx\to\Sp$ a Karoubi-localizing invariant, then, we have the following equivalence of spectra:
    $$
        F(S^1\hat\otimes\Ccal)\simeq F(\Ccal)\oplus\Sigma F(\Ccal)\oplus N_+F(\Ccal)\oplus N_-F(\Ccal)
    $$
\end{thm}

Note that we also need no longer take $\Ccal$ idempotent complete, since $\Idem(S^1\otimes\Ccal)\simeq \Idem(S^1\otimes\Idem(\Ccal))$. In particular, non-connective K-theory is Karoubi localizing and so we have the following version of the fundamental theorem for non-connective K-theory: 
$$
    \K(S^1\hat\otimes\Ccal)\simeq\K(\Ccal)\oplus\Sigma \K(\Ccal)\oplus N_+\K(\Ccal)\oplus N_-\K(\Ccal)
$$
Taking $\Ccal$ to be $\Perf(R)$ as above, where $R$ is a ring spectrum, gives a formula for non-connective K-theory of ring spectra, generalizing \cite{Grayson} and \cite{FontesOgle}. In particular, when $R=\S[\Omega X]$, this is a non-connective version of the main result of \cite{HKVWW}. \\

Taking connective covers in the formula above gives an improved version of the formula one would obtain by applying \ref{MainResultIntro} for connective K-theory, which is only Verdier-localizing. This is the actual formula appearing in \cite{Grayson} or \cite{HKVWW} in their above-mentioned specific cases. In particular, the canonical non-connective delooping of $K(\Ccal)$ appears here as the connective cover of $\Sigma\K(\Ccal)$. The connective formula given by \ref{MainResultIntro} for connective K-theory misses specifically the non-connective term appearing in $\pi_0$. \\

\textbf{Outline of the proof.} As expected, the proof of our main result relies on the idea of a \textit{projective line}, as found in \cite{HKVWW}, \cite{CisinskiKhan} or originally \cite{Grayson} but generalized to our context. There are two maps $S^1_\pm=B\N_\pm\to B\Z=S^1$ depending on the identification of $\N$ as either non-positive or non-negative integers, which induce exact functors $S^1_\pm\otimes\Ccal\to S^1\otimes\Ccal$ for any stable $\Ccal$, which we call the \textit{telescopes}. $S^1_\pm\otimes\Ccal$ models the $\Spec(A[t^{\pm1}])$ appearing in \cite{Grayson}. The pullback of the telescopes is the \textit{projective line of $\Ccal$}, denoted $\P(\Ccal)$. \\

The proof of Theorem \ref{MainResultIntro} relies on the ability to express the image of $\P(\Ccal)$ under a Verdier-localizing invariant $F$ in two different ways, one by the preservation of specific pullbacks which is a consequence of the property of being Verdier-localizing, and the other through a direct calculation, which is reliant on the fact that $\Ccal$ is idempotent complete, hence the hypothesis. This last computation is a version of the \textit{projective bundle formula} of Section 4.2 in \cite{CisinskiKhan}, whose proof has been expunged from any reference to Lurie's Spectral Algebraic Geometry (see \cite{SAG}). Indeed, we will show:

\begin{prop} \label{PropChap4Intro}
    For any Verdier-localizing invariant $F$, we have an equivalence $F(\P(\Ccal))\simeq F(\Ccal)\oplus F(\Ccal)$. Moreover, the following square is cartesian:
    $$
    \begin{tikzcd}
        F(\P(\Ccal))\arrow[r]\arrow[d] & F(S^1_+\otimes\Ccal)\arrow[d] \\
        F(S^1_-\otimes\Ccal)\arrow[r] & F(S^1\otimes\Ccal)
    \end{tikzcd}
    $$
\end{prop}

The fact that the square is cartesian is a direct, abstract consequence of the Verdier-localizing character of our functor $F$. The equivalence $F(\P(\Ccal))=F(\Ccal)\oplus F(\Ccal)$ relies on the other hand on an actual concrete calculation, going through explicit descriptions of objects at hand. This is the part that is the longer and the more intricate of the two. \\

\textbf{Organisation of the article.} The section 2 and 3 are dedicated to preliminaries regarding respectively the tensor construction and the notions of Verdier and Karoubi-localizing functors, both  outlined previously. The tensor product of section 2 is an algebraic version of the construction of 6.4.1 \cite{HermKI} and section 3 mostly extracts from the appendix of \cite{HermKII} the propositions and lemmas relevant to our problem. \\

Section 4 is where most of the magic takes place. We first define the projective line and relevant objects to prove the cartesian square part of Proposition \ref{PropChap4Intro}, and in the following subsection, we make the explicit calculation of $F(\P(\Ccal))$ for $\Ccal$ a stable $\infty$-category and $F$ Verdier-localizing. This is the most technical part of this article. \\

In section 5, we finish the proof of Theorems \ref{MainResultIntro} and \ref{MainResult2Intro}, and draw the many consequences it has for algebraic K-theory, its non-connective version as well as topological Hochschild homology and topological cyclic homology. \\

\textbf{Acknowledgements.} The author is indebted to its PhD advisor Yonatan Harpaz for suggesting the main idea of this paper, for many discussions and countless answered questions. We thank Liam Keenan and Maxime Ramzi for pointing us towards existing literature. We thank an anyonymous referee for helpful comments on an earlier draft. The author would also like to thank all of the people who have, by their unwavering and invaluable support, contributed to the completion of this article, from friends and family to colleagues and teachers. \\

Early drafts of this manuscript were written when the author was a student at the Ecole Polytechnique in Paris. The remainder was carried while the author was subsided by a PhD grant subsided by the European Research Council in the Project \textit{Foundations of motivic real K-theory}.

%%%%%%%%%%%%%%%%%%%%%%%%%%%%%%%%%%%%%%%%%%%%%%%%%%%%%%%%%%%%%%%%%%%%%%%%%%%%%%%%%%%%%%%%%%%%%%%%%%%%%%%%%%%
\section{Semi-exact tensoring of stable $\infty$-categories}

\hspace{1.2em} In this section, we are interested in a tensor-like construction between a stable $\infty$-category and a general simplicial set $K$. Ultimately, $K$ will either be $B\Z$ or $B\N$, the $\infty$-categories with one object and $\Z$ or $\N$ as (discrete) spaces of morphisms; the semi-exact tensoring of $B\N$ and $\Ccal$ will be a stable $\infty$-category which we will think as polynomials in $\Ccal$, whereas the tensoring by $B\Z$ will be Laurent polynomials. \\

Recall that a stable $\infty$-category is a pointed $\infty$-category, i.e. an $\infty$-category with a zero object, such that every morphism has a fiber and a cofiber and additionally, that fiber sequences and cofiber sequences coincide. Here, every (co)limit is to be understood in the $\infty$-categorical world and thus corresponds to a homotopy (co)limit. The study of stable $\infty$-categories is developed in the first chapter of \cite{HA}. In particular, it is shown that stable $\infty$-categories admit all of the finite limits and colimits, and cartesian squares coincide with cocartesian squares. Exact functors are functors preserving either finite limits or colimits and they in fact preserve both. \\

There are multiple constructions of $\infty$-categories which are tensor products in the correct category. For instance, when $\Ccal$ and $\Dcal$ are $\infty$-categories, $\Ccal\times\Dcal$ corepresents $\Fun(\Ccal, \Fun(\Dcal, -))$. When $\Ccal$, $\Dcal$ are stable, $\Fun^{Ex}(\Ccal, \Fun^{Ex}(\Dcal, -))$ is also corepresentable, this time in the enriched setting of $\CatEx$, and if $\Ccal$, $\Dcal$ are furthermore presentable, $\Fun^{L}(\Ccal, \Fun^{L}(\Dcal, -))$ is again internally corepresentable in $\Pr^L_{Ex}$, the $\infty$-category of presentable stable $\infty$-categories and left functors. \\

The tensor product we will use in this article is an hybridization of the first two. Indeed, we want to capture the stable nature of our stable $\infty$-categories, which we require to be able to talk about their K-theory, but we also want to tensor them by $\infty$-groupoids which need not be stable, like $B\Z$ or $B\N$. We also want the result to be stable, so that we can take its K-theory as well.

Hence, we take our tensor product to be corepresenting the functor $\Fun(K, \Fun^{Ex}(\Ccal, -))$, i.e. to be universal for functors from $K\times\Ccal$ which are exact only in the second variable.

\begin{defi}
    Let $\Ccal$ be a stable $\infty$-categories and $K$ a simplicial set. We define $K\otimes\Ccal$, the \textit{semi-exact tensoring} of $\Ccal$ by $K$, by the following universal property:
    $$
        \Fun^{Ex}(K\otimes\Ccal, \Dcal)\simeq \Fun(K, \Fun^{Ex}(\Ccal, \Dcal))
    $$
    This construction is functorial in $K$ and $\Ccal$.
\end{defi}

Note that this is \textit{not} symmetrical in $K$ and $\Ccal$ when both are stable $\infty$-categories, and does not coincide with the usual tensor product of $\CatEx$, which is universal for functors $K\times\Ccal\to\Dcal$ exact in \textit{both} variables. \\

The fact that $K\otimes\Ccal$ exists is a consequence of \cite[Proposition 5.3.6.2]{HTT}. Indeed, exact functors between stable $\infty$-categories are exactly finite-colimits preserving functors by \cite[1.1.4.1]{HA}, so $K\otimes\Ccal$ can be obtained as the universal $\infty$-category for functors out of $K\times\Ccal$ which send $\Rcal$, the finite cocones of $K\times\Ccal$ which are constant in the first variable and colimits in the second variable, to colimits at their target. This is exactly described by the construction $\Pcal^\Kcal_\Rcal(K\times\Ccal)$ given by \cite[5.3.6.2]{HTT} where $\Kcal$ designates all the finite cocones of $K\times\Ccal$. In fact, the proof of this proposition even gives an explicit description, which can be reformulated as in the following proposition:

\begin{prop} \label{TensorConstruction}
    $K\otimes\Ccal$ can be realized as the smallest subcategory of $\Fun(K^{op}, \Ind(\Ccal))$ stable by finite colimits and containing $L_{k,X}$, the left Kan extensions of $\{k\}\to \Ind(\Ccal)$ constant in $X\in\Ccal$ along the inclusions $\{k\}\subset K^{op}$. 
\end{prop}
\begin{proof}
    When $\Ccal$ has all finite colimits, $\Ind(\Ccal)$ can be identified to the full subcategory of finite colimit-preserving $\infty$-presheaves, i.e. functors $\Ccal^{op}\to S$ which preserve finite colimits. This is another consequence of \cite[5.3.6.2]{HTT}, which is explicitly stated in Example 5.3.6.8 of \textit{loc. cit}. In particular, for a stable $\Ccal$, $\Ind(\Ccal)$ is cocomplete and the left Kan extensions $L_{k,X}$ exist for any $k\in K_0$ and $X\in\Ccal$, hence our claimed construction for $K\otimes\Ccal$ is well-defined. \\
    
    Let us now unfold this construction and see why it coincides with that of \cite[5.3.6.2]{HTT}. Indeed, by what we explained above, $\Fun(K^{op}, \Ind(\Ccal))$ can be identified with functors of $\Fun(K^{op}\times\Ccal^{op}, S)$ which preserve finite colimits in the second variable. To match \cite[5.3.6.2]{HTT}, we need to show our construction identifies with the essential image of $L\circ j_0$ and is closed under finite colimits, where $L$ is the left adjoint to the inclusion $\Fun(K^{op}, \Ind(\Ccal))\subset\Fun(K^{op}\times\Ccal^{op}, S)$ and $j_0$ is the Yoneda embedding of $K\times\Ccal$. \\
    
    The closure under finite colimits is a part our definition, hence it suffices to show that $L_{k,X}\simeq L\circ j_0(k,X)$ for $k\in K_0$ and $X\in\Ccal$. Let $F$ be a functor $K\to\Ind(\Ccal)$. By \cite[4.3.3.7]{HTT}, we have the following equivalence:
    $$
        \Map_{\Fun(K^{op}, \Ind(\Ccal))}(L_{k,X}, F)\simeq \Map_{\Ind(\Ccal)}(X,F(k)) \simeq F(k)(X)
    $$
    where the second equivalence is given by the Yoneda lemma. Differently stated, this is saying that $L_{k, X}$ corepresents the evaluation functor $\Fun(K^{op}, \Ind(\Ccal))\to S$ sending $F:K^{op}\to\Ind(\Ccal)$ to $F(k)(X)$. But, by adjunction, $L\circ j_0(k,X)$ verifies:
    $$
        \Map_{\Fun(K^{op}, \Ind(\Ccal))}(L\circ j_0(k,X), F)\simeq \Map_{\Fun(K^{op}\times\Ccal^{op}, S)}(j_0(k,X), \tilde{F})
    $$
    where $j_0(k,X)$ is the Yoneda embedding at $(k, X)$ and $\tilde{F}$ is the "uncurried" functor. Hence, it follows from the Yoneda lemma that $L\circ j_0(k,X)$ is corepresenting the same functor as $L_{k,X}$, since $\tilde{F}(k, X)\simeq F(k)(X)$ by definition, and yet another instance of Yoneda gives us the wanted equivalence. This shows that $K\otimes\Ccal$ is universal for functors $K\times\Ccal\to\Dcal$ preserving finite colimits in the second variable. \\
    
    Moreover, $K\otimes\Ccal$ is a stable $\infty$-category, since it is a subcategory of the stable $\Fun(K^{op}, \Ind(\Ccal)$ which is itself stable by finite colimits and by the loop functor $\Omega$, since $\Omega L_{k, X}=L_{k, \Omega X}$ and $\Omega$ commutes with finite colimits. Since finite colimit-preserving functors between stable $\infty$-categories are exact, we have the universal property of the definition, as wanted.
\end{proof}

In \cite[Section 6.4.1]{HermKI}, the authors define the hermitian version of this construction. Because their hermitian functors are generally taken to be from $\Ccal^{op}$, their construction involves $\Pro(\Ccal)$ and right Kan extensions but this is the only difference. In particular, our proof of the proposition is in all points similar to Remark 6.4.2 establishing the universal property.

\begin{rmq} \label{TensorIsFun}
    Since $K\otimes\Ccal$ is defined by a universal property, it automatically acquires functoriality by the Yoneda lemma. Hence we have defined a functor $-\otimes-:\cat{Cat}_\infty\times\CatEx\to\CatEx$. We will give an explicit description of the induced morphisms later in the section.
\end{rmq}

\begin{ex}
    When $R$ is a ring spectrum, we will see that $S^1\otimes\Perf(R)$ identifies as a dense subcategory of $\Perf(R[t, t^{-1}])$, where $\Perf(R)$ are $R$-modules which are compact in $R\Mod$ and \textit{dense} means that every object of $\Perf(R[t, t^{-1}])$ is a retract of an object of $S^1\otimes\Perf(R)$. The same will be true for $S^1\otimes\Fun(X,\Sp)^c$ and $\Fun(X\times S^1, \Sp)^c$.
\end{ex}

In general, as the examples above show, $K\otimes\Ccal$ need not be idempotent complete even if $\Ccal$ is. However, we can identify its idempotent completion, and in fact even its $\Ind$-construction. The following lemma is a generalization of a result of \cite[Lecture 21, Proposition 6]{LurieNotes}.

\begin{lmm} \label{IdempotentCompletionTensor}
    Let $K$ be a simplicial set and $\Ccal$ an $\infty$-category, then the $\infty$-category $\Fun(K^{op}, \Ind(\Ccal))$ is compactly generated and in fact, we even have
    $$
        \Ind(K\otimes\Ccal)\simeq\Fun(K^{op}, \Ind(\Ccal))
    $$
    which means that $\Fun(K^{op}, \Ind(\Ccal))$ is generated by $K\otimes\Ccal$.
\end{lmm}
\begin{proof}
    The first claim follows from the second, since $K\otimes\Ccal$ is contained in the full subcategory of compact objects. Indeed, the left Kan extension is left adjoint to a filtered colimit-preserving functor, hence it preserves compact objects and $\Ccal$ is compact in $\Ind(\Ccal)$. Thus, it suffices to prove the announced equivalence. \\
    
    Since $\Fun(K^{op}, \Ind(\Ccal))$ is cocomplete, the inclusion $K\otimes\Ccal\subset\Fun(K^{op}, \Ind(\Ccal))$ extends to a fully-faithful $\Ind(K\otimes\Ccal)\to\Fun(K^{op}, \Ind(\Ccal))$, which we have to show is essentially surjective. But this is a map preserving colimits between presentable $\infty$-categories, hence it has a right adjoint $R$ by the adjoint functor theorem \cite[5.5.2.9]{HTT}, and it suffices to show that $R$ is conservative. \\
    
    Let $f:A\to B$ be a map of $\Fun(K^{op}, \Ind(\Ccal))$, i.e. a natural transformation between functors $K^{op}\to\Ind(\Ccal)$, such that $R(f)$ is an equivalence. Then, precomposition by $R(f)$ induces the following equivalence for any $k\in K$ and $X\in\Ccal$:
    $$
        \Map_{\Fun(K, \Ind(\Ccal))}(L_{k, X}, A)\simeq \Map_{\Fun(K, \Ind(\Ccal))}(L_{k, X}, B)
    $$
    By the universal property of left Kan extensions, it follows that
    $$
        \Map_{\Ind(\Ccal)}(X, A(k))\simeq\Map_{\Ind(\Ccal)}(X, B(k))
    $$
    for any $X\in\Ccal$. Since $\Ccal$ generates $\Ind(\Ccal)$ under filtered colimits and $X$ is compact in $\Ind(\Ccal)$, we conclude that $f$ induces an equivalence $A(k)\simeq B(k)$ for any $k\in K$. Hence $f$ is a natural equivalence as wanted.
\end{proof}

\begin{rmq} \label{IdempotentTensorProduct}
    We mentioned in introduction a third tensor-like product, for presentable stable $\infty$-categories, which is for instance the one used in \cite{Blumberg}. It induces a tensor product $\hat\otimes$ between stable $\infty$-categories $\Ccal$, $\Dcal$ which is always idempotent complete by letting $\Ccal\hat\otimes\Dcal:(\Ind(\Ccal)\otimes^L\Ind(\Dcal))^c$, where $\otimes^L$ is the symmetric monoidal structure of presentable stable $\infty$-categories with the left functors. \\
    
    We can define a hybrid version, $K\hat\otimes\Ccal:=(K\otimes^L\Ind(\Ccal))^c$ where $\otimes^L$ is defined by the following universal property\footnote{Which exists and gives a presentable stable $\infty$-category for similar reasons as in proposition \ref{TensorConstruction}}, for presentable stable $\Ccal$ and $\Dcal$:
    $$
        \Fun^L(K\otimes^L\Ccal,\Dcal)\simeq \Fun(K, \Fun^L(\Ccal,\Dcal))
    $$
    Since $\Fun^L(\Ind(\Acal), \Dcal)\simeq\Fun^{Ex}(\Acal, \Dcal)$ for stable $\Acal$ and $\Dcal$, we have that $K\otimes^L\Ind(\Ccal)=\Ind(K\otimes\Ccal)$ by comparing universal properties. Thus, we have $K\hat\otimes\Ccal \simeq \Idem(K\otimes\Ccal)$ by taking compacts objects on both sides. \\
    
    With the lemma above, we have an equivalence $K\otimes^L\Ind(\Ccal)\simeq\Fun(K^{op}, \Ind(\Ccal))$ which yields $K\hat\otimes\Ccal \simeq \Fun(K^{op}, \Ind(\Ccal))^c$. For Karoubi-localizing invariants, which we will introduce in the following section and are invariant under idempotent completion, our main result Theorem \ref{MainResultIntro} is also true for $\hat\otimes$, which is often quite easier to identify.
\end{rmq}

As we mentioned in Remark \ref{TensorIsFun}, $-\otimes\Ccal$ is functorial. Our goal now is to identify the induced exact functor $A\otimes\Ccal\to B\otimes\Ccal$ for $f:A\to B$, under the explicit description given by \ref{TensorConstruction}. 

If $f:A\to B$ is a map of simplicial set, then by \cite[4.3.3.7]{HTT} , we have an adjoint pair:
$$
    \begin{tikzcd}
        \Fun(A, \Ind(\Ccal))\arrow[r, shift left=2, "f_!"] & \Fun(B, \Ind(\Ccal))\arrow[l, shift left=2, "f^*", "\perp"']
    \end{tikzcd}
$$
where $f_!$ denotes the functor of left Kan extension along $f$ and $f*$ precomposition by $f$. Using the notations of the explicit construction of \ref{TensorConstruction}, $f_!L_{a, X}\simeq L_{f(a), X}$ for any $a\in A_0$, $X\in\Ind(\Ccal)$ since the left Kan extension of a composite is the composition of left Kan extension. Being a left adjoint, $f_!$ commutes with colimits thus restrict to an exact functor:
$$
    f\otimes\Ccal:A\otimes\Ccal\to B\otimes\Ccal
$$
which sends $L_{a, X}$ to $L_{f(a), X}$. In particular, the $\Ind$-completion of this functor is $f_!$ although remark that in general, $f^*$ need not preserve the $L_{b, X}$ so the adjunction does not descend and we only have an explicit adjoint when $\Ind$-completing. By universality, the functorial map $f\otimes\Ccal:A\otimes\Ccal\to B\otimes\Ccal$ must induce for any $\Dcal$:
$$
    \Fun^{Ex}(B\otimes\Ccal, \Dcal)\simeq \Fun(B, \Fun^{Ex}(\Ccal, \Dcal))\to \Fun(A, \Fun^{Ex}(\Ccal, \Dcal)) \simeq \Fun^{Ex}(A\otimes\Ccal, \Dcal)
$$
where the middle map is precomposition by $f$. An exact functor $F:B\otimes\Ccal\to\Dcal$ is characterized by the images $F(L_{b, X})$ for every $b\in B$ and $X\in\Ccal$, and the above map sends such an $F$ to the functor $\tilde{F}$ characterized by the data $\tilde{F}(L_{a, X})\simeq F(L_{f(a), X})$. But this is also exactly what the restriction of $f_!$ does. Hence by Yoneda, the restricted left Kan extension of $f$ is equivalent to the map $f\otimes\Ccal$ functorially induced by $f$. \\

We will need to consider some of the less well-behaved precomposition by $f_*$ in our subsequent sections. In our examples, the map $A\to B$ will be the inclusion of a unique point, and in that context, the precomposition naturally lands in $\Ind(\Ccal)$.

\begin{defi} \label{ForgetfulFunctorDef}
    Let $K$ be a 0-reduced simplicial set, then precomposition by the inclusion $*\subset K$ induces a functor $\fgt_{K}:K\otimes\Ccal\to\Fun(*,\Ind(\Ccal))\simeq\Ind(\Ccal)$ forgetting the $K$-part of the tensor. We will call it the \textit{forgetful functor} of $K\otimes\Ccal$. 
\end{defi}

\begin{rmq}
    The argument works \textit{mutatis mutandis} to show that if $F:\Ccal\to\Dcal$ is exact between stable $\infty$-categories and $K$ is a simplicial set, then the restriction of the left Kan extension functor along $F$ induces
    $$
        K\otimes F: K\otimes\Ccal\to K\otimes\Dcal
    $$
    mapping $L_{k, X}$ to $L_{k, F(X)}$, and $K\otimes F$ is indeed the map given by functoriality of $\otimes$.
\end{rmq}

%%%%%%%%%%%%%%%%%%%%%%%%%%%%%%%%%%%%%%%%%%%%%%%%%%%%%%%%%%%%%%%%%%%%%%%%%%%%%%%%%%%%%%%%%%%%%%%%%%%%%%%%%%%
\section{Verdier, Karoubi sequences and localizing invariants}

\hspace{1.2em} This section is dedicated to establishing terminology and useful results related to (split-)Verdier and Karoubi invariants. \\

These ideas were first introduced in \cite{Blumberg} under the name of \textit{additive} and \textit{localizing} invariants. However our setting fits more naturally in a middle-ground of those two notions, which was notably developed in Appendix A of \cite{HermKII}, called \textit{Verdier localizing} invariants. We adopt their terminology in the following: the \textit{localizing} invariants of \cite{Blumberg} correspond to \textit{Karoubi localizing} invariants for us, and the \textit{additive} ones to \textit{split-Verdier localizing}. Appendix A of \cite{HermKII}, which will serve as our reference of choice for this material, gives a precise comparison of all the notions in its introductory remark.

\subsection{Verdier and Karoubi sequences}

\hspace{1.2em} We recall here an array of definitions and results we will need in the following. Unless specified, they are coming from appendix A of \cite{HermKII}. First, let us define the central object of this section, \textit{Verdier sequences}.

\begin{defi}
    A \textit{Verdier sequence} is a sequence of $\CatEx$
    $$
        \begin{tikzcd}
            \Ccal\arrow[r, "i"] & \Dcal\arrow[r, "p"] &\Ecal
        \end{tikzcd}
    $$
    which is both a fiber and a cofiber in $\CatEx$. Exact functors $i$ fitting in such sequences are called Verdier inclusions and $p$ Verdier projections. \\
    
    If $p$ has a left adjoint (resp. right), then the Verdier sequence is called left-split (resp. right-split). Sequences that are both left- and right-split are called \textit{split Verdier sequences}; they are the exact sequences of \cite{Blumberg}.
\end{defi}

Recall that there is a functor $\Idem:\CatEx\to\CatExIdem$ computing the idempotent-completion of an $\infty$-category, which is left adjoint to the inclusion. In the stable setting, it preserves both limits and colimits by \cite[A.3.3]{HermKII}. Hence, sequences of $\CatEx$ that become fiber-cofibers after applying $\Idem$ form a more general class, the \textit{Karoubi sequences}:

\begin{defi}
    A \textit{Karoubi sequence} is a sequence of $\CatEx$
    $$
        \begin{tikzcd}
            \Ccal\arrow[r, "i"] & \Dcal\arrow[r, "p"] &\Ecal
        \end{tikzcd}
    $$
    which is sent by the idempotent-completion functor $\Idem$ to both a fiber and a cofiber in $\CatExIdem$. Exact functors $i$ fitting in such sequences are called Karoubi inclusions and $p$ Karoubi projections.
\end{defi}

\begin{rmq}
    Remark that the preceding definition asks the sequence to be a fiber-cofiber in $\CatExIdem$ and not $\CatEx$. Since the inclusion $\CatExIdem\subset\CatEx$ only preserves limits in general, the idempotent completion of a Verdier sequence is only a Karoubi sequence, and not a Verdier one when regarded as a sequence of $\CatEx$. In particular, there are Karoubi sequences of idempotent-complete $\infty$-categories which are \textit{not} Verdier sequences.
\end{rmq}

It will be convenient to have a way to know whether a functor $f:\Ccal\to\Dcal$ fits in a Verdier sequence. To this intent, we have the following criterion for Verdier projections and Verdier inclusions, which is extracted from A.1.6 and A.1.9 of \cite{HermKII}.

\begin{prop} \label{VerdierCriterion}
    Let $p:\Ccal\to\Dcal$ be an exact functor between stable $\infty$-categories. Then, the following are equivalent:
    \begin{enumerate}[itemsep=0em]
        \item $p$ is a Verdier projection
        \item $p$ is a localisation at some collection of arrows $\Wcal$, i.e. $p_*:\Fun(\Dcal, \Ecal)\to\Fun(\Ccal, \Ecal)$ is fully-faithful for any $\infty$-category $\Ecal$ with essential image functors $\Ccal\to\Ecal$ which invert $\Wcal$.
    \end{enumerate}
    and the following are also equivalent:
    \begin{enumerate}[itemsep=0em]
        \item $p$ is a Verdier inclusion
        \item $p$ is fully-faithful and has essential image closed under retracts in $\Dcal$
    \end{enumerate}
\end{prop}

Note that what we call localisation are the functors of Warning 5.2.7.3 of \cite{HTT} and we reserve the term \textit{Bousfield localisation} for what Lurie calls a localisation, which is asking for a fully-faithful right adjoint. 

The localizations that are Bousfield are exactly those having a right adjoint, which leads to the following criterion for left-split and right-split Verdier projections, extracted from the equivalence between $(i)$ and $(iv)$  of \cite[A.2.3]{HermKII}:

\begin{prop} \label{SplitVerdierCriterion}
    Let there be a sequence $e:\begin{tikzcd}[cramped]\Ccal\arrow[r, "f"]&\Dcal\arrow[r, "p"]&\Ecal\end{tikzcd}$ of exact functors with vanishing composite. Then, the following are equivalent:
    \begin{enumerate}[itemsep=0em]
        \item $e$ is a left-split (resp. right-split) Verdier sequence
        \item $e$ is a cofiber sequence such that $f$ is fully-faithful and has a left (resp. right) adjoint $g$
        \item $e$ is a fiber sequence such that $p$ has a fully-faithful left (resp. right) adjoint $q$
    \end{enumerate}
    When either of the propositions is satisfied, then the sequence of adjoints $\begin{tikzcd}[cramped]\Ecal\arrow[r, "q"]&\Dcal\arrow[r, "g"]&\Ccal\end{tikzcd}$ is a right-split (resp. left-split) Verdier sequence.
\end{prop}

Finally, we present the following criterion for Karoubi projections and injections, which is \cite[A.3.8]{HermKII}:

\begin{prop} \label{KaroubiCriterion}
    An exact functor $f:\Ccal\to\Dcal$ is a Karoubi injection if and only if it is fully-faithful and a Karoubi projection if it has dense image and the induced $f:\Ccal\to\im f$ is a Verdier projection.
\end{prop}

The discussion above concerns fiber-cofiber sequences and but we will need more generally properties about squares.

\begin{defi}
    A cartesian square of $\CatEx$
    $$
    \begin{tikzcd}
        \Ccal\arrow[r]\arrow[d] &\Dcal\arrow[d, "p"] \\
        \Ecal\arrow[r] & \Fcal
    \end{tikzcd}
    $$
    is said to be a split-Verdier (resp. Verdier, resp. Karoubi) square if $p$ is a split-Verdier (resp. Verdier, resp. Karoubi) projection.
\end{defi}

The following is an algebraic version of \cite[1.5.2.(iii)]{HermKII}, which is proven in the hermitian context:

\begin{lmm} \label{VerdierSquaresAreExact}
    A Verdier square is also cocartesian in $\CatEx$. The same goes for Karoubi square in $\cat{Cat}^{Idem}_\infty$ after idempotent completion.
\end{lmm}
\begin{proof}
    Given a Verdier square:
    $$
    \begin{tikzcd}
        \Ccal\arrow[r]\arrow[d, "q"] &\Dcal\arrow[d, "p"] \\
        \Ecal\arrow[r] & \Fcal
    \end{tikzcd}
    $$
    then $q$ is also a Verdier projection by \cite[A.1.11]{HermKII} , and the square extends to a diagram of cartesian squares by taking $\Gcal$ to be the fiber of $q$:
    $$
    \begin{tikzcd}
        \Gcal\arrow[r]\arrow[d] & \Ccal\arrow[r]\arrow[d, "q"] &\Dcal\arrow[d, "p"] \\
        0\arrow[r] & \Ecal\arrow[r] & \Fcal
    \end{tikzcd}
    $$
    The left-square is cartesian and cocartesian, because $q$ is a Verdier projection, and so the external rectangle is cartesian by the pasting law. Since $p$ is a Karoubi projection, the external square is also cocartesian and thus the pasting law applies to give us that the right square is cocartesian. \\
    
    The same proof works \textit{mutatis mutandis} for Karoubi squares after idempotent completion.
\end{proof}

%%%%%%%%%%%%%%%%%%%%%%%%%%%%%%%%%%%%%%%%%%%%%%%%%%%%%%%%%%%%%%%%%%%%%%%%%%%%%%%%%%%%%%%%%%%%%%%%%%%%%%%
\subsection{Verdier- and Karoubi-localizing invariants}

\hspace{1.2em} Recall that a functor $F:\Ccal\to\Dcal$ between pointed categories is called \textit{reduced} if it preserves zero objects. We are interested in the following classes of reduced functors:

\begin{defi} \label{DefVerdierLoc}
    Let $\Ecal$ be a $\infty$-category with finite limits. A reduced functor $F:\CatEx\to\Ecal$ is called:
    \begin{itemize}[itemsep=0em]
        \item \textit{additive} or \textit{split-Verdier localizing} if it sends split-Verdier squares to cartesian squares in $\Ecal$
        \item \textit{Verdier localizing} if it sends Verdier squares to cartesian squares in $\Ecal$
        \item \textit{Karoubi localizing} if it sends Karoubi squares to cartesian squares in $\Ecal$
    \end{itemize}
\end{defi}

Since split-Verdier sequences are Verdier and Verdier sequences are Karoubi, the above list is ordered so that each property implies those above it.

\begin{rmq}
    The original introduction of additive and localizing functors of \cite{Blumberg} asked furthermore that the functors preserve filtered colimits. However, this hypothesis is unused for our applications, hence we adopt a similar convention to \cite[Definition 1.5.4]{HermKII}, separating the localizing property from the filtered-colimit preservation.
\end{rmq}

When the target category $\Ecal$ is stable, our definition of split-Verdier, Verdier and Karoubi localizing invariants is actually equivalent to a weaker property, namely:

\begin{lmm}
    Let $F:\CatEx\to\Ecal$ be a reduced functor $F$ landing in a stable $\Ecal$, then $F$ is split-Verdier localizing if and only if it sends split-Verdier \textit{sequences} to exact sequences. The same applies when changing both instances of split-Verdier to Verdier or Karoubi. 
\end{lmm}
\begin{proof}
    Using the same diagram as in the proof of lemma \ref{VerdierSquaresAreExact} and applying $F$:
    $$
    \begin{tikzcd}
        F(\Gcal)\arrow[r]\arrow[d] & F(\Ccal)\arrow[r]\arrow[d] &F(\Dcal)\arrow[d] \\
        0\arrow[r] & F(\Ecal)\arrow[r] & F(\Fcal)
    \end{tikzcd}
    $$
    we see that if $F$ preserves say Verdier sequences, both the left square and the external rectangle are exact. Then by the pasting law, so is the right square, as wanted.
\end{proof}

Let $\Ccal$ be a stable $\infty$-category and denote $j:\Ccal\to\Idem(\Ccal)$ the natural map. Since $\Idem(j)$ is an equivalence, $\Ccal\to\Idem(\Ccal)\to 0$ is a Karoubi sequence. In consequence, for any Karoubi-localizing $F$, we have an equivalence $F(\Idem(\Ccal))\simeq F(\Ccal)$. In fact, the converse is true for Verdier-localizing $F$: such a $F$ is Karoubi-localizing if and only if it is invariant under idempotent completion, as we now show:

\begin{prop} \label{Karoubi=Verdier+IdempotentInvariance}
    Let $F$ be a Verdier-localizing invariant, then $F$ is Karoubi-localizing if and only if $F$ is invariant under idempotent completion.
\end{prop}
\begin{proof}
    The above discussion gives one direction of this equivalence. We adapt the proof of \cite[1.5.6]{HermKII} in our context for the other. \\
    
    Suppose $F$ Verdier-localizing and invariant under idempotent-completion. Let the following square be a Karoubi square, where $f$ and $g$ are Karoubi projections
    $$
    \begin{tikzcd}
        \Ccal\arrow[r]\arrow[d, "f"] & \Ecal\arrow[d, "g"] \\
        \Dcal\arrow[r] & \Fcal
    \end{tikzcd}
    $$
    By proposition \ref{KaroubiCriterion}, Karoubi projections factor as a Verdier projection onto its essential image followed by a fully-faithful map with dense image. Denote $\Dcal_0$ the essential image of $f$ and $\Ecal_0$ that of $g$, then the following square commutes:
    $$
    \begin{tikzcd}
        \Ccal\arrow[r]\arrow[d] & \Ecal\arrow[d] \\
        \Dcal_0\arrow[r]\arrow[d] & \Ecal_0\arrow[d] \\
        \Dcal\arrow[r] & \Fcal
    \end{tikzcd}
    $$
    where the top square is a Verdier square and the bottom vertical maps are fully-faithful with dense image. Since $F$ is invariant under idempotent completion, $\Dcal_0\to\Dcal$ is sent by $F$ to an equivalence and similarly for $\Ecal_0\to\Ecal$. But $F$ is Verdier-localizing so sends the top square to an exact square of $\Sp$. This concludes.
\end{proof}

Finally, we end this section by giving several examples of Verdier and Karoubi localizing invariants. First, we discuss algebraic K-theory and its non-connective variant. \\

Recall that algebraic K-theory $K:\CatEx\to\Sp$ and non-connective K-theory $\K:\CatEx\to\Sp$ are reduced functors preserving filtered colimits, and that the latter is invariant under idempotent completion. Moreover, if $\Ccal$ is idempotent complete, then $K(\Ccal)$ is the connective cover of $\K(\Ccal)$. We have the following:

\begin{thm} \label{KTKaroubi}
    Algebraic K-theory $K$ is Verdier localizing and non-connective K-theory $\K$ is Karoubi localizing.
\end{thm}
\begin{proof}
    By Theorem 1.3 of \cite{Blumberg}, we have that non-connective K-theory is Karoubi-localizing (since what they mean by localizing is the combination of being Karoubi-localizing and commuting with filtered colimits in our lingo), and that algebraic K-theory is split-Verdier localizing.
    
    This is enough to deduce that algebraic K-theory is Verdier localizing, as in the proof of Corollary 4.4.15 of \cite{HermKII}. We will go a slightly different route, and replace the instance of \cite{Blumberg} proving by that non-connective K-theory is Karoubi-localizing by the Special Fibration Theorem \cite[Theorem 10.20]{Barwick}, which gives a more direct proof of the fact that K-theory is Verdier-localizing. \\
    
    Suppose $\Ccal\to\Dcal\to\Ecal$ is a Verdier sequence. Then, the $\Ind$-completion of this sequence is again a Verdier sequence by \cite[A.3.11]{HermKII} (this is in fact the case for Karoubi sequences), and it is even right split by the subsequent remark of \textit{loc. cit.} with both right adjoints additionally preserving colimits. Hence, we have an accessible localisation functor $L:\Ind(\Dcal)\to\Ind(\Ecal)$ between compactly generated $\infty$-categories, induced by $\Dcal\to\Ecal$ and whose right adjoint is fully-faithful and preserves all colimits (so filtered ones in particular). Moreover, $L$-equivalences in $\Ind(\Dcal)$ are indeed generated by those between the compact objects. 
    
    Hence the Special Fibration Theorem applies, and we have a fiber sequence of spaces:
    $$
        \begin{tikzcd}
            \Kcal(\Idem(\Ccal))\arrow[r] & \Kcal(\Idem(\Dcal))\arrow[r] & \Kcal(\Idem(\Ecal))
        \end{tikzcd}
    $$
    where we identified $\Idem(C)$ with the compact objects of $\Ind(\Ccal)$ and used the letter $\Kcal$ to denote the K-theory \textit{space}. By the cofinality theorem \cite[Theorem 10.19]{Barwick}, for every stable $\Acal$, the map $K(\Acal)\to K(\Idem(\Acal))$ is injective on $\pi_0$ and an isomorphism on higher homotopy groups. In consequence, the map $K(\Ccal)\to \fib(K(\Dcal)\to K(\Ecal))$ induces an equivalence on $\pi_n$ for $n\geq 1$ by the naturality of the long exact sequence of homotopy groups. \\
    
    Since $\Dcal\to\Ecal$ is a Verdier projection, it is essentially surjective and thus $K_0(\Dcal)\to K_0(\Ecal)$ is surjective. Consequently the fiber $\fib(K(\Dcal)\to K(\Ecal))$ is connective, so it suffices to show the isomorphism of groups $K_0(\Ccal)\simeq \pi_0\fib(K(\Dcal)\to K(\Ecal))$. 
    
    The injective map $K_0(\Ccal)\to K_0(\Idem(\Ccal))$ factors through $F=\pi_0\fib(K(\Dcal)\to K(\Ecal))$ which immediately implies that $K_0(\Ccal)\to F$ is injective. In fact, the map $F\to K_0(\Idem(\Ccal))$ is also injective since it fits in the following cartesian square\footnote{The canonical map from $F$ to the pullback can be checked to be an equivalence by a diagram chase, using that maps on $K_0$ are injective and maps on $K_1$ equivalences.} of groups:
    $$
        \begin{tikzcd}
            F\arrow[d]\arrow[r] & K_0(\Dcal)\arrow[d] \\
            K_0(\Idem(\Ccal))\arrow[r] & K_0(\Idem(\Dcal))
        \end{tikzcd}
    $$
    where $K_0(\Dcal)\to K_0(\Idem(\Dcal))$ is injective. Thus, by Thomason's classification of dense subcategories \cite[Theorem A.3.2]{HermKII}, there exists a Karoubi equivalence $\widetilde{\Ccal}\to\Idem(\Ccal)$ which factors $\Ccal\to\Idem(\Ccal)$ and such that $K_0(\widetilde{\Ccal})\simeq F$.   
    
    Given the above cartesian square, we see that $\widetilde{\Ccal}$ is also the full subcategory of $X\in\Idem(\Ccal)$ such that the map $K_0(\Idem(\Ccal))\to K_0(\Idem(\Dcal))$ sends $[X]$ to the subgroup $K_0(\Dcal)$ of the target. Since $\Dcal$ is also a dense subcategory $\Idem(\Dcal)$, another instance of Thomason's theorem implies that $\widetilde{\Ccal}$ is the full subcategory of $X\in\Idem(\Ccal)$ such that $i(X)\in\Dcal$ where $i:\Idem(\Ccal)\to\Idem(\Dcal)$. But the restriction $i:\Ccal\to\Dcal$ is a Verdier inclusion hence its image is closed under retracts so $X\in\Idem(\Ccal)$ is sent to $\Dcal$ if and only if $X\in\Ccal$. This implies $\Ccal=\widetilde{\Ccal}$ which concludes. 
\end{proof}

\begin{rmq}
    A third proof of this fact has been recently given in \cite[Theorem 6.1]{HLS}. It has the notorious advantage on relying on much less machinery than the results of either Barwick or Blumberg.
\end{rmq}

We turn now our attention to topological Hochschild homology and topological cyclic homology. In \cite{Blumberg}, the authors define a functor $\THH:\CatEx\to\Sp$ which is Karoubi-localizing and commutes to filtered colimits and comes with a natural transformation $\Kth\to\THH$ called the Dennis trace map. This functors is first built on spectral categories following \cite{BlumbergMandellLocalization}, and since it inverts Morita equivalences, it induces a functor in the $\infty$-categorical setting. 

This construction on the level of spectral categories of \textit{loc. cit.} has actually more structure: it comes with a functorial $S^1$-action and even a cyclotomic structure. This refinement descends to $\infty$-categories again because it preserves Morita equivalences. Recall from Theorem II.3.7 of \cite{NikolausScholze} however that the $\infty$-categorical version of cyclotomic spectra we are considering is not quite $\CycSp$, the $\infty$-category of cyclotomic spectra of \textit{loc. cit.}, but $\CycSp^{gen}$, the $\infty$-category of genuine cyclotomic spectra, also defined in the work of Nikolaus-Scholze. Theorem II.3.8 ensures those two coincide on bounded below objects but generally we are only guaranteed an exact functor $\CycSp^{gen}\to\CycSp$. 

The upshot of this story is that we have a functor $\THH:\CatEx\to\CycSp^{gen}$, which we can further compose to get $\THH:\CatEx\to\CycSp$ and this functor recovers the $\THH$ of \cite{Blumberg} when postcomposed by the forgetful functor $\fgt:\CycSp\to\Sp$. \\

Write $\TC:\CycSp\to\Sp$ for the functor corepresented by the cyclotomic spectra $\S^{triv}$, whose underlying spectra is the sphere spectrum with trivial structure. Precomposing this functor by $\THH$ yields another functor, topological cyclic homology, which we again denote $\TC$:
$$
    \TC:\CatEx\longrightarrow\Sp
$$
By construction, $\TC$ comes with a natural transformation $\TC\to\THH$ of functors to spectra, and this natural transformation factors the Dennis trace map. The resulting factor $K\to\TC$ is called the cyclotomic trace, and features prominently in the Dundas-Goodwillie-McCarthy Theorem (see 1.2 of \cite{NikolausScholze}). 

\begin{thm} \label{THHTCKaroubi}
    Topological Hochschild homology $\THH$ and topological cyclic homology $\TC$ are Karoubi-localizing.
\end{thm}
\begin{proof}
    Proposition 10.2 of \cite{Blumberg} proves the claim for $\THH:\CatEx\to\Sp$ since their notion of localizing is stronger than our Karoubi-localizing one. Given the aforementioned refinement to cyclotomic spectra and since $\CycSp\to\Sp$ is exact and detects equivalences, $\THH:\CatEx\to\CycSp$ is again Karoubi-localizing. The claim for $\TC$ directly follows since $\TC:\CycSp\to\Sp$ is exact.
    
    Another proof of the statement for $\THH:\CatEx\to\Sp$ can be assembled from Theorem 3.4 and Proposition 4.24 of \cite{HSS}, which identify the trace of the endofunctor $\id:\Ind(\Ccal)\to\Ind(\Ccal)$ in $\Pr^L_{Ex}$ equipped with the Lurie tensor-product with the spectrum $\THH(\Ccal)$.
\end{proof}

%%%%%%%%%%%%%%%%%%%%%%%%%%%%%%%%%%%%%%%%%%%%%%%%%%%%%%%%%%%%%%%%%%%%%%%%%%%%%%%%%%%%%%%%%%%%%%%%%%%%%%%%%%%
\section{The Projective Line}
\subsection{The Projective Line as a pullback}

\hspace{1.2em} Let $S^1$ be the $\infty$-groupoid corresponding to the 1-sphere, i.e. the coherent nerve of $\Z$ seen as a discrete simplicial category. We let $S^1_+$ and $S^1_-$ be the $\infty$-categories corresponding to the inclusions of monoids $\N_+\to\Z$ and $\N_-\to\Z$. These are the same $\infty$-categories\footnote{They are $\infty$-categories because they are fibrant simplicial $\infty$-categories, since discrete simplicial spaces are Kan complexes} and only their identification within $S^1$ differs. \\

\begin{defi}
    By functoriality, the functors $S^1_\pm\to S^1$ induce exact functors $T_\pm:S^1_\pm\otimes\Ccal\to S^1\otimes\Ccal$, which we call the \textit{telescopes}. 
\end{defi}

\begin{rmq} \label{TelescopeFilteredColimit}
    Realizing the tensor as the explicit construction given by Proposition \ref{TensorConstruction} identifies the telescope $T_+$ as the left Kan extension along $S^1_+\to S^1$ of functors $S^1_+\to\Ind(\Ccal)$. Since $S^1_+\simeq B\N$ and $S^1\simeq B\Z$, this left Kan extension can be explicitly constructed by freely inverting the action of $t$ the generator of $\N$. This can be done via the standard procedure for freely adding inverses, namely taking the filtered colimit in $\Ind(\Ccal)$ of the repeated action of $t$ on $V(*)$ for a functor $V:S^1_+\to\Ind(\Ccal)$, and endowing it with the natural action of $t$ which is now invertible.
\end{rmq}

Pulling back along the two telescopes $T_\pm:S^1_\pm\otimes\Ccal\to S^1\otimes\Ccal$ yields a stable $\infty$-category that we call the \textit{Projective Line}.

\begin{defi}
    Let $\P(\Ccal)$ the $\infty$-category defined by the following pullback square:
    $$
        \begin{tikzcd}
            \P(\Ccal)\arrow[r]\arrow[d] & S^1_-\otimes\Ccal\arrow[d] \\
            S^1_+\otimes\Ccal\arrow[r] & S^1\otimes\Ccal
        \end{tikzcd}
    $$
    The $\infty$-category $\P(\Ccal)$ is stable by \cite{HA} 1.1.4.2 and we call it the \textit{Projective Line} of $\Ccal$.
\end{defi}

\begin{rmq} \label{ObjectsProjectiveLine}
    Since the inclusion $\CatEx\to\cat{Cat}_\infty$ preserves limits, we can see $\P(\Ccal)$ as a pullback in $\cat{Cat}_\infty$, where pullbacks enjoy an explicit description as homotopy limits for the Joyal model structure on $\cat{sSet}$, the category of simplicial set. Hence, $\P(\Ccal)$ has objects triples $(Y_-, Y, Y_+)$ with given equivalences $T_\pm(Y_\pm)\simeq Y$, where $Y_+\in S^1_+\otimes\Ccal$, $Y_-\in S^1_-\otimes\Ccal$ and $Y\in S^1\otimes\Ccal$. 
    % This data is also given with the necessary coherence information ..?
\end{rmq}

We now show algebraic K-theory, and in fact more generally, any Verdier-localizing invariants sends the square defining the Projective line to a cartesian square of spectra. By the theory of section 3, it suffices to show the square is a Verdier square, i.e. that the telescopes are Verdier projections.

\begin{lmm} \label{TelescopeVerdier}
    The telescopes $T_\pm:S^1_{\pm}\otimes\Ccal\to S^1\otimes\Ccal$ are Verdier projections.
\end{lmm}
\begin{proof}
    We show more generally that if $K\to L$ is a localisation functor, then $K\otimes\Ccal\to L\otimes\Ccal$ is a Verdier projection. This is in particular the case for the telescopes, since they are induced by $B\N\simeq S^1_{\pm}\to S^1\simeq B\Z$. This is a version of \cite{HermKII} 1.4.10 (i) in the algebraic context. \\
    
    By the universal property of tensor categories, we have the following commutative square for every stable $\Dcal$:
    $$
        \begin{tikzcd}
            \Fun^{Ex}(L\otimes\Ccal, \Dcal)\arrow[r, "\sim"]\arrow[d] & \Fun(L, \Fun^{Ex}(\Ccal, \Dcal))\arrow[d] \\
            \Fun^{Ex}(K\otimes\Ccal, \Dcal)\arrow[r, "\sim"]& \Fun(K, \Fun^{Ex}(\Ccal, \Dcal))
        \end{tikzcd}
    $$
    When $K\to L$ is a localization inverting a class of arrows $\Wcal$, the right vertical map is fully-faithful with essential image functors $K\to\Fun^{Ex}(\Ccal, \Dcal)$ sending $\Wcal$ to natural equivalences. Since horizontal maps are equivalences, the left vertical map is also fully-faithful and its essential image is exactly functors $K\otimes\Ccal\to\Dcal$ which send arrows of $\Wcal'$ to equivalences, where $\Wcal'$ is the class of arrows induced by a pair $(f, \id)$ with $f\in\Wcal$ and $\id$ the identity of some object in $\Ccal$. In consequence, $K\otimes\Ccal\to L\otimes\Ccal$ is a localization. By proposition \ref{VerdierCriterion}, this concludes.
\end{proof}

\begin{cor} \label{ProjectiveLineCartesianSquare}
    The cartesian square defining $\P(\Ccal)$ is a Verdier square. In particular, for any Verdier localizing invariant $F$, the following square of spectra is cartesian:
    $$
        \begin{tikzcd}
            F(\P(\Ccal))\arrow[r]\arrow[d] & F(S^1_-\otimes\Ccal)\arrow[d] \\
            F(S^1_+\otimes\Ccal)\arrow[r] & F(S^1\otimes\Ccal)
        \end{tikzcd}
    $$
\end{cor}

%%%%%%%%%%%%%%%%%%%%%%%%%%%%%%%%%%%%%%%%%%%%%%%%%%%%%%%%%%%%%%%%%%%%%%%
\subsection{An explicit calculation for the Projective Line}

\hspace{1.2em} The preceding section gave a straightforward understanding of $\P(\Ccal)$ under a Verdier localizing invariant since such a pullback square is preserved by such an invariant. The goal of this section is to show it is possible to make an effective calculation of $F(\P(\Ccal))$ under slightly stricter hypotheses on either our Verdier invariant $F$ or our stable $\infty$-category $\Ccal$. Let us first give a definition

\begin{defi} \label{DefPsiZero}
    Let $\Ccal$ be a stable category. There is an exact functor $\psi_0:\Ccal\to\P(\Ccal)$ which is given by the left Kan extensions $L_{*, X}^\pm$ and $L_{*, X }$ on the respective components. It is well-defined because we have equivalences $T_\pm(L_{*, X}^\pm)\simeq L_{*, X }$ which are natural in $X$.
\end{defi}

Our computation of $F(\P(\Ccal))$ is given by the following theorem:

\begin{thm} \label{ProjectiveLineEffectiveCalculation}
    Let $F:\CatEx\to\Sp$ be a Verdier-localizing invariant and $\Ccal$ a stable idempotent-complete $\infty$-category. The natural map $\psi_0:\Ccal\to\P(\Ccal)$ is a right-split Verdier inclusion with cofiber $\Ccal$, which in particular induces an equivalence of spectra $F(\P(\Ccal))\simeq F(\Ccal)\oplus F(\Ccal)$.
\end{thm}

\begin{rmq}
    The above theorem is our version of the \textit{projective bundle formula} of Theorem 4.2.5 in \cite{CisinskiKhan}. Our proof will be significantly longer as we cannot (and do not want to) rely on the spectral algebraic geometry developed in \cite{SAG}. 
\end{rmq}

The remainder of the section is dedicated to the proof. Fix a stable $\infty$-category $\Ccal$, which we will suppose furthermore idempotent complete after lemma \ref{PhiPreservesC}. To construct the equivalence of the theorem, we introduce a functor $\Phi:\P(\Ccal)\to\Ccal$ which we will show is a left-split Verdier projection. We will then identify its fiber as none other than $\Ccal$, and the equivalence of the theorem will follow from the splitting lemma. \\

For the construction of $\Phi$, recall the explicit description of the objects of $\P(\Ccal)$ of \ref{ObjectsProjectiveLine}. Objects are triplets $(Y_-, Y, Y_+)$ coming with equivalences $T_\pm(Y_\pm)\simeq Y$. Equivalently, by the adjunction between left Kan extensions and precompositions, these are maps $Y_\pm\to i_\pm(Y)$ of $\Fun(S^1_\pm,\Ind(\Ccal))$, where $i_\pm$ is the forgetful functor associated to $S^1_\pm\to S^1$ (see Definition \ref{ForgetfulFunctorDef}). \\

In consequence, to each object of $\P(\Ccal)$ is functorially associated a map $Y_-\oplus Y_+\to Y$ in $\Ind(\Ccal)$ where $Y_-$, $Y$ and $Y_+$ should be their image by the respective forgetful functors, which we abusively denoted the same way to avoid unnecessary notation clutter. This defines a functor $\P(\Ccal)\to\Map(\Delta^1, \Ind(\Ccal))$.

\begin{defi}
    Let $\Phi:\P(\Ccal)\to \Ind(\Ccal)$ be the functor defined as the composite 
    $$
        \begin{tikzcd}
        \P(\Ccal)\arrow[r] & \Fun(\Delta^1, \Ind(\Ccal))\arrow[r, "\fib"] & \Ind(\Ccal)
        \end{tikzcd}
    $$ 
    On objects, this is the fiber of $Y_-\oplus Y_+\to Y$, namely, we have the following exact square in $\Ind(\Ccal)$:
    $$
        \begin{tikzcd}
            \Phi(Y_-, Y, Y_+)\arrow[r]\arrow[d] & Y_-\oplus Y_+\arrow[d] \\
            0\arrow[r] & Y
        \end{tikzcd}
    $$
\end{defi}

In \cite{HKVWW}, whose proof serves as inspiration for ours, the authors use a slightly different construction: instead of the functor $\Phi$, they introduce a functor $\Gamma$, the global section functor, obtained by taking the cofiber of the map $Y_-\oplus Y_+\to Y$ instead of the fiber. This was necessary in the context of \textit{loc. cit.} since the Waldhausen categories considered were not supposed to have all finite limits. In our stable setting where both exist, we found the fiber to be easier to work with and hence replaced instances of $\Gamma$ by $\Phi$. This is merely by convenience: all of the following could be done by replacing $\Phi$ by $\Gamma$ and suitably changing the proofs. \\

We also need to consider shift functors, which are the higher categorical version of those in \cite{HKVWW}:
\begin{defi} \label{DefShifts}
    Let $n\in\Z$, the \textit{$n$-shift functor} $[n]:\P(\Ccal)\to\P(\Ccal)$ is the functor given on objects by sending a triple $(Y_-, Y, Y_+)$ to the same triple but where the equivalence $T_-(Y_-)\simeq Y$ is composed by the equivalence $t^n:Y\to Y$.
\end{defi}

There is an arbitrary choice made in working with the $S^1_-$-component for shifts, and one could define a shift on the $S^1_+$-side. However, a triple $(Y_-, Y, Y_+)$ shifted on the $Y_+$-side is equivalent in $\P(\Ccal)$ to a shift on the $Y_-$-side of the same triple, the equivalence being induced by the following commutative diagram:
$$
    \begin{tikzcd}[row sep=large, column sep=large]
        T_-(Y_-)\arrow[r]\arrow[d, equal] &  Y\arrow[d, "t^{-n}"] & T_+(Y_+)\arrow[l, "t^n"]\arrow[d, equal] \\
        T_-(Y_-)\arrow[r, "t^{-n}"] &  Y & T_+(Y_+)\arrow[l]
    \end{tikzcd}
$$
In the following, whenever we mention a shift on the $S^1_+$-components, we also implicitly apply the above equivalence to get a $S^1_-$-shift. \\

Using shifts, one can lift any map between the $S^1\otimes\Ccal$-components of objects of $\P(\Ccal)$ to a map of $\P(\Ccal)$ itself.

\begin{lmm} \label{LiftMaptoProjective}
    Let $y=(Y_-, Y, Y_+)$ and $z=(Z_-, Z, Z_+)$ be two objects of $\P(\Ccal)$ and $F:Y\to Z$. Then, there exists $z'\in\P(\Ccal)$ which only differs from $z$ by shifts and a map $f:y\to z'$ whose component $Y\to Z$ is the prescribed $F$ up to a shift of a power of $t$. 
\end{lmm}
\begin{proof}
    We appeal to the description of the telescope of Remark \ref{TelescopeFilteredColimit}. After forgetting the extra structure, i.e. as an object of $\Fun((S^1_-)^{op}, \Ind(\Ccal))$, $Z$ is the filtered colimit of $Z_-$ under the action of $t$. But $Y_-$ is compact in $\Fun((S^1_-)^{op}, \Ind(\Ccal))$ since it belongs to $S^1_-\otimes\Ccal$, hence the map $Y_-\to Z$ induced by the equivalence $T_-(Y_-)\simeq Y$ factors through $Z_-$ at some finite point, i.e. through some finite shift of $Z_-$. This, with the analogous procedure for $Y_+$ and $Z_+$ (see the remark above), gives a globally finite object $z'$ with the same components as $z$ but shifted equivalences $T_\pm(Z_\pm)\simeq Z$, and $z'$ comes with a well-defined map $f:y\to z'$ with the wanted $F:Y\to Z$ up to a shift by a power of $t$ as its middle component. 
\end{proof}

The crucial observation is that $\Phi$ lands in the subcategory of $\Ind(\Ccal)$ of compact objects, namely $\Idem(\Ccal)$. In particular, when $\Ccal$ is idempotent complete, this is the functor $\P(\Ccal)\to\Ccal$ we were looking for.

\begin{lmm} \label{PhiPreservesC}
    For any $y\in\P(\Ccal)$, we have $\Phi(y)\in\Idem(\Ccal)$.
\end{lmm}
\begin{proof}
    Any object $X\in\Ccal$ gives rises to an object of $\P(\Ccal)$ by taking respective left Kan extensions $L_{*, X}^\pm$ and $L_{*, X}$ of the constant functor $*\to\Ind(\Ccal)$ associated to $X$ along the inclusions of the point in $S^1_\pm$ and respectively $S^1$. Those left Kan extensions are compatible with the telescopes, i.e. there are equivalences $T_\pm(L^\pm_{*, X})\simeq L_{*, X}$ in $S^1\otimes\Ccal$. Hence, the triple $\psi_0(X):=(L_{*, X}^+, L_{*, X}, L_{*, X}^-)$ actually defines an object of $\P(\Ccal)$.

    The $\infty$-category of such objects and their shifts need unfortunately not be stable, but we can consider $\P^{gf}(\Ccal)$, the smallest full stable subcategory of $\P(\Ccal)$ containing the above objects as well as their shifts. We call objects of $\P^{gf}(\Ccal)$ \textit{globally finite}. \\

    Our first claim is that $\Phi$ maps globally finite objects to $\Ccal$. Indeed, once we have forgotten the extra structure, then we simply have $L_{*, X}\simeq\coprod_{\Z}X$ and $L_{*, X}^\pm\simeq\coprod_{\N_\pm}X$ as objects of $\Ind(\Ccal)$, and the map $L_{*, X}^+\oplus L_{*, X}^-\simeq L_{*, X}$ is
    $$
        \coprod_{\N_+}X\oplus\coprod_{\N_-}X\longrightarrow\coprod_{\Z}X
    $$
    which is the inclusion of each summand along the inclusions $\N_\pm\subset\Z$. Thus, upon shifting once, we have an equivalence\footnote{We remind here the reader that what we call $\N$ is a monoid, hence contains 0, in good Bourbaki fashion, and thus there is a double identification of the zeroth summand in the above map.} which means $\Phi((L_{*, X}^+, L_{*, X}, L_{*, X}^-)[-1])\simeq 0$ and shifting either adds copies of $X$ or copies of $\Omega X$, depending on whether copies of $X$ are missed or mapped-to twice. Since $\Phi$ is an exact functor, it preserves finite colimits and thus sends globally finite objects to $\Ccal$. \\
    
    Let $y=(Y_-,Y, Y_+)\in\P(\Ccal)$. We claim we can find a globally finite object $z=(Z_-, Z, Z_+)$ with a map $F=(f_-, f, f_+):y\to z$ such that the component $f:Y\to Z$ is an equivalence. Indeed, by lemma \ref{LiftMaptoProjective}, it suffices to build a globally finite object $z$ with an equivalence $F:Y\to Z$. Letting $Z=Y$ and $F$ be the identity, this reduces to building a globally finite object whose middle component is any prescribed $Z\in S^1\otimes\Ccal$. 
    
    Using that $S^1\otimes\Ccal$ is the smallest full stable subcategory of $\Fun(S^1, \Ind(\Ccal))$ containing the free $L_{*, X}$ for $X\in\Ccal$ and that $(L_{*, X}^+, L_{*, X}, L_{*, X}^-)$ is a suitable globally finite object with middle term $L_{*, X}$, the above property follows from the stability by pushout of "being the middle component of a globally finite object". Given a pushout square
    $$
        \begin{tikzcd}
            X\arrow[r]\arrow[d] & Y\arrow[d] \\
            Z\arrow[r] & T
        \end{tikzcd}
    $$
    and globally finite objects $(X_-, X, X_+)$, $(Y_-, Y, Y_+)$ and $(Z_-, Z, Z_+)$, one can lift the span $Z\longleftarrow X\longrightarrow Y$ to $\P^{gf}(\Ccal)$ using Lemma \ref{LiftMaptoProjective} up to shifting the given globally finite objects. Then, taking the actual pushout in $\P^{gf}(\Ccal)$ is done component-wise (the inclusion $\P^{gf}(\Ccal)\subset\P(\Ccal)$ is exact by definition), hence gives a globally finite object $(T_-, T, T_+)$ whose middle component is the prescribed $T$. This proves our claim. \\
    
    Having chosen such a globally finite object $z$ with a map $F$, we can consider the cokernel of $F$, which is an object of $\P(\Ccal)$ obtained as the following triple $\coker(F)\simeq(\coker(f_-), 0, \coker(f_+))$. In consequence, it verifies
    $$
        \Phi(\coker(F))=\coker(f_-)\oplus\coker(f_+)
    $$
    where the $\coker(f_\pm)$, originally computed in $S^1_\pm\otimes\Ccal$, are now considered as objects of $\Ind(\Ccal)$ (i.e. there is an implicit forgetful functor $\fgt_{S^1_\pm}$). \\
    
    We now show that in fact, $\coker(f_\pm)\in\Idem(\Ccal)$. In $S^1_\pm\otimes\Ccal$, we have a natural map $\coker(f_-)\to T_-(\coker(f_-))$ induced by the identity of $T_-(\coker(f_-)$. Since $T_-(\coker(f_-))=0$ and $\coker(f_-)$ is compact in $\Fun(S^1_-, \Ind(\Ccal))$, the filtered colimit description guarantees that there is some $m\geq 0$ such that $t^m:\coker(f_-)\to\coker(f_-)$ is the zero map. In consequence, the identity $\coker(f_-)\to\coker(f_-)$ factors through $\coker(t^m)$. But $\coker(t^m)\in\Ccal$, indeed this is clear if $\coker(f_-)$ is of the form $L_{*, X}^-$ and in general, $\coker(f_-)$ is a finite colimit of such objects, which commute to the formation of $\coker(t^m)$. This shows the $\coker(f_\pm)$ are retracts of objects of $\Ccal$, hence are in $\Idem(\Ccal)$. In consequence, $\Phi(\coker(F))\in\Idem(\Ccal)$ since the latter is stable by direct sum. \\
    
    Since $\Phi$ preserves finite colimits, the following square is exact:
    $$
        \begin{tikzcd}
            \Phi(y)\arrow[r, "\Phi(F)"]\arrow[d] & \Phi(z)\arrow[d] \\
            0\arrow[r] & \Phi(\coker(F))
        \end{tikzcd}
    $$
    We have $\Phi(z)\in\Ccal$ and $\Phi(\coker(F))\in\Idem(\Ccal)$ hence the fiber $\Phi(y)\in\Idem(\Ccal)$ as wanted. This concludes.
\end{proof}

Hence, we have a well-defined map $\Phi:\P(\Ccal)\to\Idem(\Ccal)$. In particular, when $\Ccal$ is idempotent complete, \textbf{which we now suppose for the rest of the section}, we have a functor $\Phi:\P(\Ccal)\to\Ccal$. Our goal is first to show that $\Phi$ is a left-split Verdier projection and then to identify its fiber with $\Ccal$.

We will show $\psi_0$, which we defined in Definition \ref{DefPsiZero}, is the wanted left adjoint of $\Phi$. As a preliminary remark, we have the following calculation:

\begin{lmm} \label{CalculationPhiPsi0}
    We have an equivalence $\Phi\circ\psi_0(X)\simeq X$ for any $X\in\Ccal$.
\end{lmm}
\begin{proof}
    By definition, $\Phi(\psi_0(X))$ is the fiber of $\coprod_{\N_+}X\oplus \coprod_{\N_-}X\to \coprod_{\Z}X$. This is simply the projection $X\oplus \coprod_{\Z}X\to \coprod_{\Z}X$, which has fiber $X$.
\end{proof}

We can now show the following:

\begin{prop}
    $\Phi$ is a left-split Verdier projection.
\end{prop}
\begin{proof}
    In order to show this, we only need showing $\Phi$ has a fully-faithful left adjoint, by \ref{SplitVerdierCriterion}. We show it is in fact given by the functor $\psi_0:\Ccal\to\P(\Ccal)$. Since $\Phi\circ\psi_0\simeq \id$ by Lemma \ref{CalculationPhiPsi0}, it will follow from the adjunction that $\psi_0$ is fully-faithful. \\
    
    The functor $\Phi$ is actually the restriction of a more general $\Psi:\Ind(\P(\Ccal))\to\Ind(\Ccal)$, where $\Ind(\P(\Ccal))$ is the following pullback:
    $$
    \begin{tikzcd}
        \Ind(\P(\Ccal))\arrow[r]\arrow[d] & \Fun(S^1_+, \Ind(\Ccal))\arrow[d] \\
        \Fun(S^1_-, \Ind(\Ccal))\arrow[r] & \Fun(S^1, \Ind(\Ccal))
    \end{tikzcd}
    $$
    and $\Psi$ is given by the same formula. Clearly, $\Psi$ has a left adjoint $\Psi_0:\Ind(\Ccal)\to\Ind(\P(\Ccal))$ which gives component-wise the left Kan extension of a point $X\in\Ind(\Ccal)$ along the inclusion $*\to S^1_\pm$.
    
    Restricting $\Psi$ to $\P(\Ccal)$ gives $\Phi:\P(\Ccal)\to\Ccal$ since $\Ccal$ is idempotent-complete and restricting $\Psi_0$ to $\Ccal$ yields $\psi_0:\Ccal\to\P(\Ccal)$ hence the adjunction descends to the restrictions.
    %%% TODO: if unclear, revive this:
    %We can decompose $\Phi$ as the following composite
    %$$
    %\begin{tikzcd}
    %    \P(\Ccal)\arrow[r] & \Ind(\Ccal)\oplus\Ind(\Ccal)\oplus\Ind(\Ccal)\arrow[r] & \Fun(\Delta^1, \Ind(\Ccal))\arrow[r, "\fib"] & \Ind(\Ccal)
    %\end{tikzcd}
    %$$
    %The second functor has a left adjoint given by the constant functor, so it remains to investigate the first functor. The latter is itself a composite of the component-wise forgetful functor to $\Ccal$, which are each right adjoints to free functors (left Kan extension from a point), and then taking the direct sum, which is right adjoint once again to the constant functor.
\end{proof}

Since $\Phi$ is a left-split Verdier projection, the following fiber sequence is a left-split Verdier sequence:
$$
    \begin{tikzcd}
        \P(\Ccal)^\Phi\arrow[r] & \P(\Ccal)\arrow[r, "\Phi"] & \Ccal
    \end{tikzcd}
$$
where the superscript $\Phi$ indicates the fiber of $\Phi$. Hence, this sequence will split under a Verdier localizing invariant. Thus, to prove the theorem, it now suffices to identify $\P(\Ccal)^\Phi$ with $\Ccal$. \\

Recall the definitions of the $n$-shifts functors $[n]$ given in \ref{DefShifts}. Clearly, $[n]$ and $[-n]$ are inverses of one another, and the proof of the lemma \ref{CalculationPhiPsi0} shows $\psi_0(X)[-1]$ lies in the fiber of $\Phi$. Thus we have an adjoint pair $\Phi\circ[1]$ and $[-1]\circ\psi_0$ which descends to $\P(\Ccal)^\Phi$ and $\Ccal$. Since $[-1]\circ\psi_0$ is again fully-faithful, this means $\Phi\circ[1]$ is also a right-split Verdier projection and we have the following right-split Verdier sequence:
$$
    \begin{tikzcd}
        (\P(\Ccal)^\Phi)^{\Phi\circ[1]}\arrow[r] & \P(\Ccal)^\Phi\arrow[r] & \Ccal
    \end{tikzcd}
$$
where the superscript $\Phi\circ[1]$ denotes again the fiber. This right-split Verdier sequence will once again split under a Verdier localizing invariant. Thus to conclude, it suffices to show the fiber $(\P(\Ccal)^\Phi)^{\Phi\circ[1]}$ is zero, which is done in the following lemma:

\begin{lmm} \label{DoubleFiberLemma}
    For any $\Ccal$ stable, we have $(\P(\Ccal)^\Phi)^{\Phi\circ[1]}\simeq 0$.
\end{lmm}
\begin{proof}
    Let $y=(Y_-, Y, Y_+)\in\P(\Ccal)$, there is a map $u:y\to [1]y$ which is the identity on $Y$ and $Y_+$ and $t^{-1}:Y_-\to Y_-$ on the first term. There is also another map $d:y\to [1]y$ which has the identity on $Y_-$ and $t:Y\to Y$ as well as $t:Y_+\to Y_+$. Those maps can be shifted by $[-1]$ and we claim the composites $d\circ (u[-1])$ and $u\circ (d[-1]):[-1]y\to [1]y$ are equivalent; indeed, this is canonically the case on non-negative and non-positive parts as it amounts to a choice of homotopy $f\circ\id\simeq\id\circ f$ for $f$ the action of either $t$ or $t^{-1}$, consequently the two homotopies also agree on the total space hence providing a global homotopy between the two maps in $\P(\Ccal)$. Hence we have a commutative square
    $$
        \begin{tikzcd}
            {}[-1](y)\arrow[r, "{u[-1]}"]\arrow[d, "{d[-1]}"'] & {}(y)\arrow[d, "u"] \\
            {}(y)\arrow[r, "d"'] & {}[1](y)
        \end{tikzcd}
    $$
    This square is in fact exact; by \cite{HTT} 5.4.5.5, it suffices that the two projections on $S^1_\pm\otimes\Ccal$ of this square are exact, which is clear since horizontal morphisms are identities and vertical morphisms identical. \\
    
    The functor $[n]$ preserves cartesian squares since it is an equivalence, hence we more generally have a cartesian square for every $n\in\N$:
    $$
        \begin{tikzcd}
            {}[n-1](y)\arrow[r, "{u[n-1]}"]\arrow[d, "{d[n-1]}"'] & {}[n](y)\arrow[d, "{u[n]}"] \\
            {}[n](y)\arrow[r, "{d[n]}"'] & {}[n+1](y)
        \end{tikzcd}
    $$
    Since $\Phi$ is exact, it preserves pullbacks and applying it to the above square, we have a third exact square:
    $$
        \begin{tikzcd}
            \Phi([n-1](y))\arrow[r]\arrow[d] & \Phi([n](y))\arrow[d] \\
            \Phi([n](y))\arrow[r] & \Phi([n+1](y))
        \end{tikzcd}
    $$
    Thus, if $y\in\P(\Ccal)$ is such that both $\Phi(y)\simeq 0$ and $\Phi\circ[1](y)\simeq 0$, i.e. an element of $(\P(\Ccal)^{\Phi})^{\Phi\circ[1]}$, then induction on $n\in\Z$ implies that $\Phi\circ[n](y)\simeq 0$ for every $n$. Hence, when $y\in(\P(\Ccal)^\Phi)^{\Phi\circ[1]}$, the following exact sequence
    $$
        \begin{tikzcd}[column sep=large]
            \Phi\circ[n](y)\arrow[r] & Y_-\oplus Y_+\arrow[r, "t^n\alpha_-\oplus\alpha_+"] & Y
        \end{tikzcd}
    $$
    gives that $t^n\alpha_-\oplus\alpha_+$ is an equivalence. \\
    
    To conclude, our proof now takes a detour through homotopy groups of a stable $\infty$-category $\Ccal$, so let us recall quickly what they are. For $X\in\Ccal$, we denote $\pi_{X}(-):=\pi_0\Map_{\Ind(\Ccal)}(X,-)$, the zeroth homotopy group of the spectrum\footnote{Recall that stable $\infty$-categories are naturally enriched in $\Sp$} $\Map_{\Ind(\Ccal)}(X, -)$. For $n\in\Z$, the stability of $\Ccal$ implies that taking $\pi_n$ instead of $\pi_0$ in the preceding formula simply changes $\pi_{X}$ to $\pi_{\Sigma^nX}$. In consequence, since $\Ind(\Ccal)$ is generated under filtered colimits by $\Ccal$, the Yoneda lemma implies that the $\pi_{X}$ jointly detect equivalences for $X\in\Ccal$. \\
    
    Using remark \ref{TelescopeFilteredColimit}, the equivalences $T_\pm(Y_\pm)\simeq Y$ show that $Y$ is the filtered colimit of a tower of $Y_\pm$ where the maps are induced by the action of $t^{\pm 1}$. For any $X\in\Ccal$, $\Map_{\Ind(\Ccal)}(X, -)$ preserves filtered colimits because $X$ is compact in $\Ind(\Ccal)$, hence for any $v_+\in\pi_{X}(Y_+)$, its image $\alpha_+(v_+)$ in $\pi_X(Y)$ can be realised as some $t^n\alpha_-(v_-)$ for some $v_-\in\pi_X(Y_-)$ and some $n\in\N_+$. But then, $\alpha_+\oplus t^n\alpha_-$ sends $(v_+, -v_-)$ to zero in $\pi_X(Y)$. Since it is also an isomorphism by above, this means $v_+$ was zero to start with. \\
    
    Hence, $\pi_{X}(Y_+)=0$ for all compact $X$, meaning $Y_+\simeq0$ and dually $Y_-\simeq0$. In consequence, $Y\simeq T_+(0)\simeq 0$, and finally $y=0$. Hence we have shown that $(\P(\Ccal)^\Phi)^{\Phi\circ[1]}$ is zero as wanted.
\end{proof}

This concludes the proof of Theorem \ref{ProjectiveLineEffectiveCalculation}.

%%%%%%%%%%%%%%%%%%%%%%%%%%%%%%%%%%%%%%%%%%%%%%%%%%%%%%%%%%%%%%%%%%%%%%%%%%%%%%%%%%%%%%%%%%%%%%%%%%%%%%%%%%
\section{The Fundamental Theorem of Verdier-localizing invariants and its consequences}
\vspace{0.5em}

\hspace{1.2em} Assembling the results from the preceding section gives us the following: for any Verdier-localizing $F$ and any stable idempotent-complete $\Ccal$, the following square is cartesian:
$$
    \begin{tikzcd}
        F(\Ccal)\oplus F(\Ccal)\arrow[r]\arrow[d] & F(S^1_-\otimes\Ccal)\arrow[d] \\
        F(S^1_+\otimes\Ccal)\arrow[r] & F(S^1\otimes\Ccal)
    \end{tikzcd}
$$
Our first subsection shows how to turn this cartesian square into the announced Theorem \ref{MainResultIntro} (\ref{MainResult} in the text) and its corollary Theorem \ref{MainResult2Intro} when $F$ is furthermore Karoubi-localizing (\ref{MainResult2} in the text). In the subsequent subsections, we draw consequences from this Theorem for algebraic K-theory of spaces, algebraic K-theory of rings and topological Hochschild homology as well as topological cyclic homology.

%%%%%%%%%%%%%%%%%%%%%%%%%%%%%%%%%%%%%%%%%%%%%%%%%%%%%%%%
\subsection{Main results}

\hspace{1.2em} Combining the results of the previous section, we have the following:
\begin{thm} \label{MainResult}
    Let $\Ccal$ be a stable idempotent complete $\infty$-category and $F:\CatEx\to\Sp$ a Verdier-localizing invariant. Then, we have the following equivalence of spectra:
    $$
        F(S^1\otimes\Ccal)\simeq F(\Ccal)\oplus\Sigma F(\Ccal)\oplus N_+F(\Ccal)\oplus N_-F(\Ccal)
    $$
    where $N_\pm F(\Ccal)$ are equivalent nil-terms.  
\end{thm}
\begin{proof}
    From the preceding section, we have a cartesian square
    $$
        \begin{tikzcd}
            F(\Ccal)\oplus F(\Ccal)\arrow[r]\arrow[d] & F(S^1_-\otimes\Ccal)\arrow[d] \\
            F(S^1_+\otimes\Ccal)\arrow[r] & F(S^1\otimes\Ccal)
        \end{tikzcd}
    $$
    The top left corner is obtained by composing the square of \ref{ProjectiveLineCartesianSquare} and the equivalence $F(\Ccal)\oplus F(\Ccal)\simeq F(\P(\Ccal))$ of Theorem \ref{ProjectiveLineEffectiveCalculation} induced by $\psi_0$. Remark that by Lemma \ref{CalculationPhiPsi0}, both arrows $F(\Ccal)\to F(S^1_\pm\otimes\Ccal)$ have a retraction given by the map induced by $\Phi$, hence they split in $\Sp$. We thus have an equivalence:
    $$
        F(S^1_\pm\otimes\Ccal)\simeq F(\Ccal)\oplus N_\pm F(\Ccal)
    $$
    where $N_\pm F(\Ccal)$ is the fiber of the respective splitting map. For the same reasons, $F(S^1\otimes\Ccal)$ splits as $F(\Ccal)\oplus\Pcal$ with some fiber $\Pcal$. The maps $F(\Ccal)\oplus F(\Ccal)\to F(\Ccal)\oplus N_\pm F(\Ccal)$ are by definition  zero on the nil-term. Hence, taking the fiber by the first projection, we have the following cartesian square:
    $$
        \begin{tikzcd}
            F(\Ccal)\arrow[r]\arrow[d] & N_-F(\Ccal)\arrow[d] \\
            N_+F(\Ccal)\arrow[r] & \Pcal
        \end{tikzcd}
    $$
    where both maps $F(\Ccal)\to N_\pm F(\Ccal)$ are zero. All of the above construction are natural hence the above square defines an exact sequence of functors $F\to N_+F\oplus N_-\to \Pcal$, where the first map is nullhomotopic. This means we have in fact a commutative diagram in $\Fun(\CatEx,\Sp)$:
    $$
        \begin{tikzcd}
            F\arrow[r]\arrow[d] & 0\arrow[r]\arrow[d] & N_+F\oplus N_-F\arrow[d] \\
            0\arrow[r] & \Sigma F\arrow[r] & \Pcal
        \end{tikzcd}
    $$
    The left square is exact by definition and the outer one by the above. The pasting law then implies the right square is also exact. Since the suspension is computed pointwise, we have the following equivalence, natural in stable idempotent-complete $\Ccal$:
    $$
        \Pcal\simeq\Sigma F(\Ccal)\oplus N_-F(\Ccal)\oplus N_+F(\Ccal)
    $$
    Given that $F(S^1\otimes\Ccal)\simeq F(\Ccal)\oplus\Pcal$, this concludes the proof of the main theorem. 
\end{proof}

\begin{rmq}
    By Lemma \ref{TelescopeVerdier}, the map $S^1_+\otimes\Ccal\to S^1\otimes\Ccal$ is a Verdier projection, hence if we denote $\Nil(\Ccal)$ its fiber in $\CatEx$, we have for every Verdier localizing $F:\CatEx\to\Sp$
    $$
        F(\Nil(\Ccal))\simeq F(\Ccal)\oplus\Omega N_+F(\Ccal)
    $$
    A similar splitting appeared in the setting of additive 1-category with twisted coefficients in \cite{LuckSteimle}.
\end{rmq}

In the specific case where $F$ is Karoubi-localizing, which is exactly asking it is Verdier-localizing and invariant under idempotent completion by Proposition \ref{Karoubi=Verdier+IdempotentInvariance}, our statement works more generally for any stable $\Ccal$. This is because $S^1\otimes\Ccal$ and $S^1\otimes\Idem(\Ccal)$ have the same idempotent completion, namely $S^1\hat\otimes\Ccal\simeq\Fun(S^1, \Ind(\Ccal))^c$. Hence in that case, we can further replace $\otimes$ by $\hat\otimes$ in the formula of \ref{MainResult}.

\begin{thm} \label{MainResult2}
    Let $\Ccal$ be \textit{any} stable $\infty$-category and $F:\CatEx\to\Sp$ a Karoubi-localizing invariant. Then, we have the following equivalence of spectra:
    $$
        F(S^1\otimes\Ccal)\simeq F(\Ccal)\oplus\Sigma F(\Ccal)\oplus N_+F(\Ccal)\oplus N_-F(\Ccal)
    $$
    where $N_\pm F(\Ccal)$ are equivalent nil-terms extended from Theorem \ref{MainResult} in the obvious way. \\
    
    Since $\Idem(S^1\otimes\Ccal)=S^1\hat\otimes\Ccal=\Fun(S^1, \Ind(\Ccal))^c$, we also have the following equivalence of spectra:
    $$
        F(S^1\hat\otimes\Ccal)\simeq F(\Ccal)\oplus\Sigma F(\Ccal)\oplus N_+F(\Ccal)\oplus N_-F(\Ccal)
    $$
    which is often the more practical formula of the two.
\end{thm}

In the rest of this section, we will exclusively draw consequences from Theorem \ref{MainResult2}. As a first corollary, recall that Theorem \ref{KTKaroubi} states that non-connective K-theory $\K$ is Karoubi localizing, hence the following formula:

\begin{cor} \label{CorollaryKTheory}
    We have the following equivalence of spectra for a stable $\Ccal$:
    $$
        \K(S^1\hat\otimes\Ccal)\simeq\K(\Ccal)\oplus\Sigma\K(\Ccal)\oplus N_+\K(\Ccal)\oplus N_-\K(\Ccal)
    $$
\end{cor}
    
\begin{rmq}
    Theorem \ref{KTKaroubi} also states that (connective) algebraic K-theory is Verdier-localizing, so one may be tempted to apply Theorem \ref{MainResult} directly, without going through the non-connective version. However, the formula obtained this way differs in $\pi_0$ from the usual formula of the fundamental theorem of K-theory. \\
    
    Indeed, the formula misses the term induced by the non-connective delooping, since the suspension of a connective spectra is always connected. This does not make the formula incorrect (thankfully!), however, since $K(S^1\otimes\Ccal)$ is not $K(S^1\hat\otimes\Ccal)$, and differs exactly by a factor in $\pi_0$ per Wall's finiteness obstruction (see \cite{LurieNotes} lecture 15, Theorem 17 for a modern version of this). \\
    
   This exemplifies the reason we will mostly be using this second version in the following, even when we want results about connective K-theory which is not Karoubi-localizing: the usual formulas of the literature such as \cite{Grayson} or \cite{HKVWW} involve non-connective terms which cannot be obtained without the non-connective input of negative K-groups, i.e. non-connective K-theory. Still, the first formula can find some use, notably in establishing a version of the formula for the finite version of algebraic K-theory of space.
\end{rmq}

%%%%%%%%%%%%%%%%%%%%%%%%%%%%%%%%%%%%%%%%%%%%%%%%%%%%%%%%
\subsection{The Fundamental Theorem for algebraic K-theory of spaces}

\hspace{1.2em} In this short subsection, we explain how Corollary \ref{CorollaryKTheory} can be used to extend the fundamental theorem of algebraic K-theory of spaces, proved in \cite{HKVWW}, to a non-connective version. From this, we will also be able to deduce their version for the connective K-theory of spaces, and actually extend it to the context of spectra. \\

We shortly recall how algebraic K-theory of spaces is defined in our $\infty$-categorical context (see for instance Lecture 21 of \cite{LurieNotes}). Let $X$ be a space, which in the following will mean a simplicial set which is an $\infty$-groupoid. We are interested in two stable $\infty$-categories $\Fun(X, \Sp)^{fin}$ and its idempotent completion $\Fun(X, \Sp)^{c}$. The latter is the subcategory of compact objects of $\Fun(X, \Sp)$ and the former is the smallest stable full subcategory containing all the left Kan extension along the inclusion $*\subset BX$ of constant functors to $\Sp$.

We denote $A^{f}(X)$ the connective K-theory of the first and $A^{fd}(X)$ that of the second; they are respectively the \textit{finite} and \textit{finitely dominated} version of the algebraic K-theory of the space $X$. Taking non-connective algebraic K-theory of either yields a non-connective K-theory of the space $X$, which we denote $\A^{fd}(X)$. The original version of the A-functor, defined by Waldhausen in \cite{Waldhausen}, is the finite version and the finitely dominated one, which only differs in its $\pi_0$ (a phenomenon related to Wall's finiteness obstruction, see \cite{LurieNotes}), was introduced in \cite{HKVWW} in order to prove the same formula we are now going to produce. \\

By lemma \ref{IdempotentCompletionTensor} and the subsequent remark, if $X$ is a space, then the idempotent completion of $S^1\otimes\Fun(X, \Sp)^c$ is $\Fun(X\times S^1, \Sp)^c$ since $\Ind(\Fun(X, \Sp)^c)\simeq \Fun(X, \Sp)$. Moreover, the explicit construction of the tensor product given by Proposition \ref{TensorConstruction} gives an equivalence $S^1\otimes\Fun(X, \Sp)^{fin}\simeq \Fun(X\times S^1, \Sp)^{fin}$. \\

Applying Corollary \ref{CorollaryKTheory} with $\Ccal=\Fun(X, \Sp)^c$ then gives the following:

\begin{thm} \label{FundTHM-ATheory}
    Let $X$ be a space. Then, we have the following splitting of non-connective K-theory:
    $$
        \A^{fd}(S^1\times X)\simeq\A^{fd}(X)\oplus\Sigma\A^{fd}(X)\oplus N_+\A^{fd}(X)\oplus N_-\A^{fd}(X)
    $$
    where $N_\pm\A^{fd}(X)$ are equivalent nil-terms.
\end{thm}

Taking the connective cover of the theorem yields the known formula of \cite{HKVWW} since all categories in question are idempotent complete. Namely, we have:

\begin{cor} \label{FundTHM-ATheoryHKVWW}
    Let $X$ be a space. Then, we have the following splitting of finitely-dominated algebraic K-theory of spaces:
    $$
        A^{fd}(S^1\times X)\simeq A^{fd}(X)\oplus BA^{fd}(X)\oplus N_+A^{fd}(X)\oplus N_-A^{fd}(X)
    $$
    where $N_\pm A^{fd}(X)$ are equivalent nil-terms and $BA(X)$ is the (non-connective) delooping of $A(X)$ which has $\pi_{-1}\A^{fd}(X)$ as its $\pi_0$.
\end{cor}

As we explained before, it is not possible to deduce from this formula a version for finite algebraic K-theory of spaces, because we rely crucially on the idempotent-completeness of the stable $\infty$-categories in question, and the finite categories of module are not idempotent-complete.

%%%%%%%%%%%%%%%%%%%%%%%%%%%%%%%%%%%%%%%%%%%%%%%%%%%%%%%%
\subsection{The Fundamental Theorem for K-theory of ring spectra}

\hspace{1.2em} The result of the previous section falls in fact in a more general setting: non-connective K-theory of arbitrary ring spectra, i.e. $A_\infty$-ring objects of $\Sp$ aka $\E_1$-ring objects of $\Sp$. The non-connective K-theory of a stable $\infty$-category can always be expressed as the (filtered colimit of) non-connective K-theory of ring spectra, as explained in \cite{LurieNotes} (Lecture 19, Remark 4), hence this is as much a general statement as our main theorem. \\

For algebraic K-theory of a space $X$, this ring spectrum is $\S[\Omega X]$. Since $\S[\Omega X] [t, t^{-1}]\simeq \S[\Omega (X\times S^1)]$, the version of the fundamental theorem for ring spectrum, which we will now explicit, recovers as a special case the formula of the section above for algebraic K-theory of spaces. \\

We denote $R\Mod$ the $\infty$-category of $R$-module spectra, and $\Perf(R)$ its subcategory of compact objects. We let $\K(R)$ and $K(R)$ be the non-connective and regular K-theory of $\Perf(R)$. Just as for algebraic K-theory of spaces, the idempotent completion of $S^1\otimes\Perf(R)$ is $\Perf(R[t, t^{-1}])$, using that the $\Ind$-construction of $\Perf(R)$ is $R\Mod$, and denoting $R[t, t^{-1}]$ the ring spectrum of Laurent polynomials in $R$. \\

Applying Corollary \ref{CorollaryKTheory} to $\Ccal=\Perf(R)$ yields the following fundamental theorem for algebraic K-theory of ring spectra:

\begin{thm}
    For $R$ a ring spectrum, we have the following splitting of non-connective K-theory:
    $$
        \K(R[t, t^{-1}])\simeq \K(R)\oplus \Sigma\K(R)\oplus N_+\K(R)\oplus N_-\K(R)
    $$
    where $N_\pm\K(R)$ are equivalent nil-terms.
\end{thm}

In particular, taking connective covers gives a formula for the connective K-theory of a ring spectrum:

\begin{thm}
    For $R$ a ring spectrum, we have the following splitting of connective K-theory:
    $$
        K(R[t, t^{-1}])\simeq K(R)\oplus \Bcal K(R)\oplus N_+K(R)\oplus N_-K(R)
    $$
    where $N_\pm K(R)$ are nil-terms obtained as the connective covers of $N_\pm\K(R)$ and $\Bcal K(R)$ is the non-connective delooping of $K(R)$ whose $\pi_0$ is $\K_{-1}(R)$.
\end{thm}

\begin{rmq}
    If $R$ is a ring, then denote $HR$ the associated Eilenberg-MacLane ring spectrum. The K-theories $\Kth(R)$ and $\Kth(HR)$ are equivalent by Barwick's theorem of the heart \cite[Theorem 6.1]{BarwickHeart}, since $\Kth(R)$ is the K-theory of the 1-category of compact $R$-modules and $\Kth(HR)$ is the category of its bounded derived $\infty$-category which is exactly $\Perf(HR)$. \\
    
    Since $(HR)[t, t^{-1}]\simeq H(R[t, t^{-1}])$, the preceding theorem recovers the usual Bass-Heller-Swan fundamental theorem of the K-theory of rings.
\end{rmq}

Like in the case of algebraic K-theory of spaces, we could also consider the smallest stable subcategory of $R\Mod$ containing $R$, that we denote $R\Mod^{fin}$. It is contained in $\Perf(R)$, which is its idempotent completion, and taking its K-theory plays a similar role to finite algebraic K-theory of spaces with regards to its finitely-dominated variant. However, due to those categories not being idempotent-complete and connective K-theory not being Karoubi localizing, we cannot apply either of our results.

%%%%%%%%%%%%%%%%%%%%%%%%%%%%%%%%%%%%%%%%%%%%%%%%%%%%%%%%
\subsection{The fundamental theorem for $\THH$ and $\TC$}

\hspace{1.2em} By Theorem \ref{THHTCKaroubi}, topological Hochschild homology and topological cyclic homology, central tools of trace methods, are Karoubi-localizing. In consequence, Theorem \ref{MainResult2} applies to them and gives the following Bass-Heller-Swan formulas:

\begin{thm}
    Let $\Ccal$ be a stable $\infty$-category, then we have the two following splitting, the first in cyclotomic spectra and the second in spectra:
    \begin{itemize}
        \item $\THH(S^1\hat\otimes\Ccal)\simeq \THH(\Ccal)\oplus\Sigma \THH(\Ccal)\oplus N_+\THH(\Ccal)\oplus N_-\THH(\Ccal)$
        \item $\TC(S^1\hat\otimes\Ccal)\simeq \TC(\Ccal)\oplus \Sigma \TC(\Ccal)\oplus N_+\TC(\Ccal)\oplus N_-\TC(\Ccal)$
    \end{itemize}
    where $N_+\THH(\Ccal)$ and $N_-\THH(\Ccal)$ are equivalent nilterms and similarly for $N_+\TC(\Ccal)$ and $N_-\TC(\Ccal)$.
\end{thm}

In particular, all the formulas for algebraic K-theory, $\THH$ and $\TC$ are deduced from Theorem \ref{MainResult2} which is natural in the invariant, the Dennis trace map and the cyclotomic trace preserve the decompositions of the fundamental theorem. \\

The same specialisation as for algebraic K-theory above can be done for $\THH$ and $\TC$ of ring spectra. In particular, knowing that $\THH(\S)=\S^{triv}$ is the cyclotomic sphere, we find that
$$
    \THH(\S[t, t^{-1}])=\S^{triv}\oplus\Sigma\S^{triv}\oplus N_+\THH(\S)\oplus N_-\THH(\S)
$$
where the last two summands are equivalent. As a cyclotomic spectrum, the nilterm $N\THH(\S)$ is of importance: indeed, using \cite{NikolausScholze}, we see that for a cyclotomic spectrum $X$, we have an equivalence $TR(X)\simeq\Map_{\CycSp}(N\THH(\S), X)$. The underlying spectrum of $N\THH(\S)$ is known to split as an infinite coproduct (see \cite{McCandless} where it is called $\tilde{\THH}(\S[t])$).

\bibliographystyle{alpha}
\bibliography{bibliographie}

\end{document}